\documentclass{amsart}
\input epsf

\usepackage{graphicx}

\usepackage{fullpage}

\usepackage{amsfonts,amsthm,amsmath,amssymb,latexsym}

\def\RR{\mathbb{R}}

\def\T{\mathcal{T}}

\def\S{\mathcal{H}}
\def\tmu{\tilde{\mu}}

\begin{document}
\title{Teichm\"uller theory of the punctured solenoid}
\author{R. C. Penner and Dragomir \v{S}ari\'{c}}

\address{Departments of Mathematics and Physics/Astronomy,
University of Southern California, Los Angeles, CA 90089}
\email{rpenner@math.usc.edu}

\address{Institute for Mathematical Sciences, Stony Brook University,
Stony Brook, NY 11794-3660}
\email{saric@math.sunysb.edu}


\date{\today}

\begin{abstract}
The punctured solenoid $\S$ is an initial object for the category
of punctured surfaces with morphisms given by finite covers
branched only over the punctures. The (decorated) Teichm\"uller
space of $\S$ is introduced, studied, and found to be parametrized
by certain coordinates on a fixed triangulation of $\S$.
Furthermore, a point in the decorated Teichm\"uller space induces
a polygonal decomposition of $\S$ giving a combinatorial
description of its decorated Teichm\"uller space itself. This is
used to obtain a non-trivial set of generators of the modular
group of $\S$, which is presumably the main result of this paper.
Moreover, each word in these generators admits a normal form, and the
natural equivalence relation on
normal forms is described. There is furthermore a
non-degenerate modular group invariant two form on the
Teichm\"uller space of $\S$.  All of this structure is in perfect analogy with that
of the decorated Teichm\"uller space of a punctured surface of finite type.
\end{abstract}

\maketitle

%VT
\thispagestyle{empty}
\def\IMSmarkvadjust{0 pt}
\def\IMSmarkhadjust{0 pt}
\def\IMSmarkhpadding{0 pt}
\def\IMSpubltext{Published in modified form:}
\def\SBIMSMark#1#2#3{
 \font\SBF=cmss10 at 10 true pt
 \font\SBI=cmssi10 at 10 true pt
 \setbox0=\hbox{\SBF \hbox to \IMSmarkhpadding{\relax}
                Stony Brook IMS Preprint \##1}
 \setbox2=\hbox to \wd0{\hfil \SBI #2}
 \setbox4=\hbox to \wd0{\hfil \SBI #3}
 \setbox6=\hbox to \wd0{\hss
             \vbox{\hsize=\wd0 \parskip=0pt \baselineskip=10 true pt
                   \copy0 \break%
                   \copy2 \break% 
                   \copy4 \break}}
 \dimen0=\ht6   \advance\dimen0 by \vsize \advance\dimen0 by 8 true pt
                \advance\dimen0 by -\pagetotal
	        \advance\dimen0 by \IMSmarkvadjust
 \dimen2=\hsize \advance\dimen2 by .25 true in
	        \advance\dimen2 by \IMSmarkhadjust

%
%   Check for publication info
%
%  \newread\jref
  \openin2=publishd.tex
  \ifeof2\setbox0=\hbox to 0pt{}
  \else 
     \setbox0=\hbox to 3.1 true in{
                \vbox to \ht6{\hsize=3 true in \parskip=0pt  \noindent  
                {\SBI \IMSpubltext}\hfil\break
                \input publishd.tex 
                \vfill}}
  \fi
  \closein2
  \ht0=0pt \dp0=0pt
 \ht6=0pt \dp6=0pt
 \setbox8=\vbox to \dimen0{\vfill \hbox to \dimen2{\copy0 \hss \copy6}}
 \ht8=0pt \dp8=0pt \wd8=0pt
 \copy8
 \message{*** Stony Brook IMS Preprint #1, #2. #3 ***}
}

\SBIMSMark{2005/06}{August 2005}{}

\vskip 1cm
\section{Introduction}

Sullivan \cite{Sul} introduced the {\it universal hyperbolic solenoid}
as the inverse limit of all finite unbranched covers of a
compact surface with negative Euler characteristic. (Different
choices of compact surface give homeomorphic solenoids since any two such
surfaces admit a common cover.) The space of
all complex structures on the solenoid is a version of a
``universal'' Teichm\"uller space insofar as the union of Teichm\"uller
spaces of all compact surfaces lies naturally as a dense subset
\cite{Sul},
\cite{NS}.

Locally as a topological space, the solenoid is modeled on a
product of the form (2-dimensional disk)~$\times$~(Cantor set),
and these charts glue together to provide a 2-dimensional
foliation of the solenoid itself. Because surface groups are
residually finite, each leaf of this foliation is a 2-dimensional
disk, which is in fact dense in the solenoid. Deformations of
geometric structures on the solenoid are typically required to be
smooth/conformal/quasiconformal in the 2-disk direction and
continuous in the Cantor set direction.  One may follow the
pattern of Ahlfors-Bers theory \cite{Ah1}, \cite{AB} for compact
surfaces in order to precisely define the Teichm\"uller space and
modular group of the solenoid \cite{Sul}.

The Ehrenpreis Conjecture, which is well-known in certain circles, is
that any two compact Riemann surfaces have almost
conformal finite unbranched covers of the same genus.  Sullivan
\cite{Sul}, \cite{BN} noted that the Ehrenpreis Conjecture is equivalent
to the statement that the (baseleaf preserving) modular group of the
solenoid
has dense orbits in the Teichm\"uller
space of the solenoid.  The algebraic structure of the modular
group of Sullivan's solenoid is not yet well understood. (An interesting
phenomenon is that any finite subgroup of the modular group of the
solenoid is cyclic
\cite{MS}, unlike for compact surfaces.)

We modify the universal object by allowing controlled finite branching,
namely, the {\it punctured solenoid} $\S$ (our universal object) is the
inverse limit of all finite unbranched covers of any fixed punctured
surface with negative Euler characteristic, e.g., covers of the
once-punctured torus, where, in effect, branching is permitted only over
the missing puncture. Equivalently again using properties of finite
covers, $\S$ is the inverse limit over all finite-index subgroups
$K$ in the modular group $PSL_2({\mathbb Z})$ of the tower ${\mathbb D}/K$
of covers, where
${\mathbb D}$ is the unit disk with frontier circle $S^1$.

Unlike Sullivan's universal hyperbolic solenoid, $\S$ is not a compact
space; the ends are quotients of the product (horoball)$\times$(Cantor
set) by the continuous action of a countable group, and the orbit of each
horoball is dense in the end.  The centers of the horoballs are called
punctures of $\S$.

In analogy to the case of punctured surfaces \cite{P2}, we
introduce decorations at the punctures, namely, a choice of
horocycle at each puncture, and we find global coordinates and a
combinatorial decomposition of the decorated Teichm\"uller space
of $\S$.  Furthermore, the combinatorial action of an appropriate
``baseleaf preserving'' modular group of the punctured solenoid is
used to give an explicit set of generators. The elements of this modular group
are written as compositions of generators in ``normal form''. We find the relations which identify
equivalent normal forms. We furthermore disprove the Ehrenpreis
Conjecture in a strong sense for the decorated Teichm\"uller
space, namely, there is an open dense subspace of the quotient by
the modular group which is Hausdorff.  Finally, we give a closed two form akin to the Weil-Petersson
K\"ahler form using our coordinates, and we show it is invariant under the modular group using our generators.

With this statement of the main results of this paper complete, we turn to
a somewhat more detailed description and discussion.

The Teichm\"uller space $T(\S )$ of the punctured solenoid $\S$ is
a separable complex Banach space \cite{Sul} with a complete
Teichm\"uller metric which is equal to the induced (by the complex
structure) Kobayashi metric \cite{Sa}. Starting from the quasiconformal definition
of Teichm\"uller space, we use the universal
covering space and covering group for $\S$ from
\cite{Sa} to give a
representation-theoretic definition of
$T(\S )$ in the spirit of hyperbolic geometry (see Theorem 4.1).

The decorated
Teichm\"uller space
$\tilde{T}(\S )$ is the space of decorated hyperbolic structures on $\S$,
i.e., all hyperbolic metrics on $\S$ together with an assignment of one
horocycle centered at each puncture. The analogous decorated Teichm\"uller
spaces of punctured surfaces as well as a universal space
modeled on orientation-preserving homeomorphisms of the circle modulo the
M\"obius group were studied
\cite{P2}, \cite{P1} using certain coordinates adapted to the
decorated setting called {\it lambda lengths}; the lambda
length is essentially the hyperbolic distance between horocycles
(precisely the square root of twice the exponential of this signed
distance, taken with a positive sign when the horocycles are disjoint).

For a
punctured surface, the decorated Teichm\"uller space is parametrized by
all positive assignments of lambda lengths, one such coordinate for each
edge of any fixed ``ideal triangulation'' of the surface, i.e., a
decomposition into triangles with vertices at the punctures (cf.
Theorem B.1).

In the universal setting, the corresponding coordinates are
assigned to the edges of the {``Farey tesselation''} of ${\mathbb
D}$, which we next recall and for which we establish our notation.
Fix the ideal hyperbolic triangle $T$ with vertices $1$, $-1$ and
$-i$ on the boundary $S^1$ of ${\mathbb D}$ sitting in the complex
plane, and let $e_0$ denote the real segment connecting $\pm 1$.
The group generated by hyperbolic reflections in the sides of $T$
contains $PSL_2({\mathbb Z})$ as the subgroup of
orientation-preserving elements.  The orbit of $T$ under this
group gives the edges of the Farey tesselation of ${\mathbb D}$ as
illustrated in Figure~1.  $\tau _*$ is regarded as a set of edges.
The ideal vertices of $\tau _*$, i.e., the asymptotic points in
$S^1$, are naturally enumerated by $\bar{\mathbb Q}={\mathbb
Q}\cup\{\infty\}\subseteq S^1$ as is also illustrated. In fact,
$PSL_2({\mathbb Z})$ acts simply transitively on the set of
orientations on edges of $\tau _*$. The edges of $\tau _*$ are
also naturally enumerated by $\mathbb{Q}-\{ 1,-1\}$: any edge
$e\in\tau _*-\{ e_0\}$ separates two triangles, whose vertices are
given by the endpoints of $e$ and two other points in $\bar
{\mathbb Q}$; among the two latter, take as the label for $e$ the
one on the opposite side of $e$ from the edge $e_0$ of $\tau _*$.
We shall also typically regard $e_0$ as a distinguished oriented
edge, or ``{\it DOE}'', from $-1$ to $+1$ as is also illustrated
in Figure~1.

\hskip .9in
\includegraphics{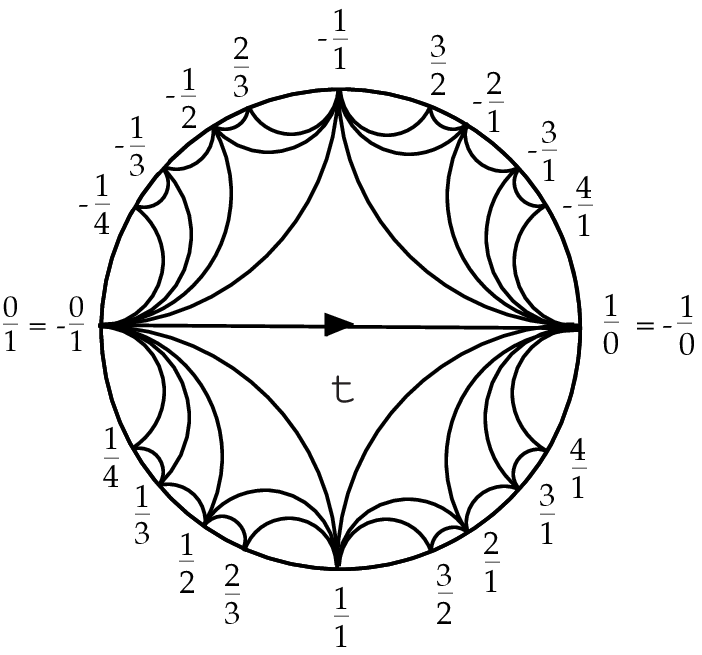}

%\epsffile{farey.eps}

\centerline{{\bf Figure~1}~The Farey tesselation $\tau _*$.}

\vskip .3cm

In the universal setting of \cite{P1}, if $f$ is a homeomorphism
of the circle, then we derive another ``tesselation with {\it
DOE}'' $\tau _f$ of ${\mathbb D}$ by demanding that three points
in $S^1$ span an ideal triangle complementary to $\tau _f$ if and
only if their pre-images under $f$ span a triangle complementary
to $\tau _*$.  Using rigidity of tesselations with {\it DOE} and
properties of order-preserving maps of the circle, it is not hard
to see that every tesselation (i.e., locally finite decomposition
by geodesics into ideal triangles) arises in this way.    This
essentially proves that the space of orientation-preserving
homeomorphisms of the circle is identified (after the choice of
$\tau _*$ as a kind of basepoint) with the collection of all
tesselations with {\it DOE} of ${\mathbb D}$.  We may then
decorate the ideal points of $\tau _f$ (i.e., choose one horocycle
centered at each point of $\tau _f$) and introduce one lambda
length coordinate to each edge of $\tau _*$.

In contrast to the case of
punctured surfaces,
an
assignment of putative positive lambda lengths to each edge of $\tau _*$
does not necessarily correspond to a point in the universal decorated
Teichm\"uller space, and it is an open question
which weights do
\cite[Section 3]{P1}; on the other hand, a sufficient condition from
\cite[Theorems 6.3 and 6.4]{P1} is that if the lambda lengths are
``pinched'' in the sense that there is a constant
$M>1$ with all lambda lengths between
$M^{-1}$ and $M$, then there is a corresponding homeomorphism of the
circle, and in fact, this homeomorphism is quasi-symmetric (the
latter being joint with Sullivan in
\cite{P1}). This is related to an open question by Thurston
\cite{Th}, namely, which measured laminations of the unit disk produce an
earthquake whose boundary values give a homeomorphism of the unit circle.

Letting $G$ denote the fundamental group of the punctured surface used
to define the punctured solenoid and letting $\hat G$ denote its
profinite completion (see Section 2), the set $Cont^G(\hat{G}, {{\mathbb R}}_{>0}^{\tau
_*})$ of all continuous
$G$-equivariant weights $\hat G\to {\mathbb R}_{>0}^{\tau _*}$
(see  Section 5 for the precise definition)
is in bijection with
the decorated Teichm\"uller space $\tilde{T}(\S )$, and the natural norm
on
$Cont^G(\hat{G}, {{\mathbb R}}_{>0}^{\tau _*})$ renders this map a
homeomorphism. These are important properties of
$\tilde{T}(\S )$ in view of the problematic behavior of other universal
Teichm\"uller spaces in this regard.

\vskip .2 cm

\noindent{\bf Theorem 5.3}~~\it The assignment of lambda lengths
is a homeomorphism onto
$$\lambda :\tilde{T}(\S_G )\to Cont^G(\hat{G}, {{\mathbb R}}_{>0}^{\tau
_*});$$ that is, we obtain a parametrization of
$\tilde{T}(\S_G)$.\rm

\vskip .2 cm

Whereas in the universal setting of \cite{P1}, there are lambda length
numbers assigned to each edge of $\tau _*$, for the solenoid, there are
continuous lambda length functions $\hat G\to{\mathbb R}_{>0}$ for each
edge of $\tau _*$.  Furthermore, whereas the weak topology on the former
corresponds to the compact-open topology on the space of homeomorphisms of
the circle, it is the strong topology on lambda length functions which
corresponds to the decorated Teichm\"uller theory of the solenoid.

The ``convex hull construction'' \cite{P2}  is the basic construction in the
decorated setting which provides combinatorial from geometric data.
In the case of punctured surfaces, it gives a modular group invariant
cell decomposition of decorated Teichm\"uller space, where cells in the
decomposition are in one-to-one correspondence with decompositions
of the surface into ideal polygons, or ``pavings'' (cf. Appendix B).  In
the universal setting of
\cite{P1}, the pinched condition on lambda lengths was shown to
be sufficient to guarantee that the corresponding decomposition of
${\mathbb D}$ is again a locally finite one by ideal
polygons with a possibly infinite number of sides
(cf. Appendix A), and we shall again call such a decomposition a
``paving'' of
${\mathbb D}$. It turns out that continuity of the lambda length functions
easily implies that the lambda lengths are pinched for each fixed element
of $\hat G$ (see Lemma 6.1), so
these aspects of
\cite{P1} directly apply to the punctured solenoid.

We may choose a leaf of the foliation of $\S$ designated the
``baseleaf''.  If
$\tau$ is a  paving of ${\mathbb D}$, let
$\stackrel{\circ}{\mathcal C}(\tau )$ be the set of decorated
hyperbolic structures on $\S$ so that the convex hull construction
associates the decomposition $\tau$ of ${\mathbb D}$ on the baseleaf.  As
$\tau$ varies, this gives a decomposition of $\tilde{T}(\S)$ which is
invariant under the modular group.  In contrast to the case of punctured
surfaces (cf. Theorem~B.5), it is not known whether each decomposition
element for the solenoid is contractible, and indeed, it is not even known
precisely which  are non-empty.

A tractable class of pavings of ${\mathbb D}$ is provided by the
``transversely locally constant'', or TLC, ones.  Namely, choose a
subgroup $K$ of $PSL_2({\mathbb Z})$ of finite index without
elliptics, choose a paving of the surface ${\mathbb D}/K$,
and lift to the universal cover to get a $K$-invariant paving $\tau$ of
${\mathbb D}$ which is said to be TLC.  These are obviously a very
special class of pavings of ${\mathbb D}$, and yet they are generic in
the following sense:

\vskip .2 cm

\paragraph{\bf Theorem 6.2}{\it The subspace $\stackrel{\circ}
{\mathcal C}(\tau )$ of $\tilde T({\mathcal H})$ is open for each
TLC triangulation $\tau$, and
 $\cup _\tau~\stackrel{\circ} {\mathcal
C}(\tau )$ is a dense open subset of $\tilde T({\mathcal H})$,
where the union is over all TLC triangulations $\tau$ of $\S$.}

\vskip .2 cm

The {\it baseleaf preserving modular group} $Mod_{BLP}(\S )$
consists of all quasiconformal baseleaf preserving self-maps of
$\S$ up to isotopies. One achievement of this paper is to give an
explicit generating set for $Mod_{BLP}(\S )$ together with
an explicit ``normal form'' for  elements of $Mod_{BLP}(\S )$.
To this end, we first show that $Mod_{BLP}(\S )$ acts transitively
on $\{\stackrel{\circ} {\mathcal C}(\tau ),\ \tau\mbox{ is a TLC
triangulation of }\S\}$ and find the isotropy groups.

\vskip .2 cm

\paragraph{\bf Theorem 7.6} {\it $Mod_{BLP}(\S_G )$ acts
transitively on $\{ \stackrel{\circ}{\mathcal C}(\tau):\tau ~{\rm
is~a~TLC}~{\rm tesselation}\}$.  Furthermore, the isotropy subgroup in
$Mod_{BLP}(\S_G )$ of $\stackrel{\circ}{\mathcal C}(\tau )$ is
quasi-conformally conjugate to $PSL_2({\mathbb Z})$.}

\vskip .2 cm

Fix an ideal triangulation $\tau$ of a punctured surface $S$. The
Whitehead move on an edge $e\in\tau$ replaces $e$ by the other
diagonal of the quadrilateral in $(S-\tau )\cup {e}$ and keeps the
remaining edges of $\tau$ fixed. In \cite{P2} (cf. Theorem B.4), a
finite presentation was given for $Mod(S)$
generated by Whitehead
moves and symmetries of top-dimensional cells, where certain
sequences of Whitehead moves correspond to elements of $Mod(S)$; however,
a single Whitehead move typically does not correspond to an element of
$Mod(S)$.

We appropriately define Whitehead moves on any TLC triangulation
$\tau$ of $\S$ in effect by performing the move $K$-equivariantly for
some finite-index subgroup $K$ as before.  In this context, a single
Whitehead move is  an element of
the modular group.  A sequence of these elements is said to be
``geometric'' if there is actually a sequence of Whitehead moves
along consecutive ideal triangulations underlying it.  (See Section 8.)
Moreover, we show that geometric sequences of Whitehead moves and
$PSL_2(\mathbb{Z})$ generate the modular group of $\S$. A non-trivial
generating set of the modular group for Sullivan's universal hyperbolic
solenoid is not known.

\vskip .2 cm

\noindent {\bf Theorem 8.5} {\it Any element of the modular group
$Mod_{BLP}(\S_G )$ can be written as a composition $w\circ
\gamma$, where $\gamma\in PSL_2({\mathbb Z})$ and $w$ is a
geometric composition of $K$-equivariant Whitehead homeomorphisms
for some fixed $K$.  }

\vskip .2 cm

We distinguish four relations on the generators above. For the
detailed description, see Section 8. Three of the relations are
already familiar from \cite{P1},\cite{P2}: the pentagon,
commutativity and involutivity relations. (The fourth relation is related to cosets
and is specific to the punctured solenoid.)  In
fact, the composition $\omega\circ\gamma$ in the above theorem is
called a {\it normal form} of the element of $Mod_{BLP}(\S )$, and we
next give a necessary condition for equivalent normal forms.

\vskip .2 cm

\noindent {\bf Theorem 8.7}\it ~~Suppose $\gamma_1 ,\gamma_2\in
PSL_2(\mathbb{Z} )$ and  $\omega_1 ,\omega_2$ are two geometric
Whitehead compositions with $\omega_1\circ\gamma_1
=\omega_2\circ\gamma_2$.  Then there is some finite index subgroup
$K$ of $PSL_2({\mathbb Z})$ with $S={\mathbb D}/K$ so that up to
Relations 1-4) for $K$, $\omega _2^{-1}\circ\omega _1=\gamma
_2^{-1}\circ \gamma _1$ is a finite composition $\phi _1\circ\phi
_2\circ\cdots\circ \phi _k$ of automorphisms $\phi _i\in Aut(\tau
_i)<PSL_2({\mathbb Z})$ of tesselations $\tau _i$ of $S$ without
{\it DOE}, for $i=1,\ldots, k$. \rm

\vskip .2 cm

We show that an orbit of a point under $Mod_{BLP}(\S )$ is not
dense in $\tilde{T}(\S )$, namely the Ehrenpreis Conjecture is
false for $\tilde{T}(\S )$. More precisely, we show that the space $\cup
_\tau~\stackrel{\circ} {\mathcal C}(\tau )/Mod_{BLP}(\S_G )$ is Hausdorff,
where the union is over all TLC triangulations $\tau$.

\vskip .2 cm

\paragraph{\bf Theorem 7.7} {\it The quotient $\cup_{\tau}
\stackrel{\circ}{\mathcal C}(\tau ) /Mod_P(\S )$ is Hausdorff,
where the union is over all TLC triangulations $\tau$. Moreover,
no orbit of a point in $\tilde{T}(\S )$ is dense.}

\vskip .2 cm

We finally introduce a generalization of the Weil-Petersson K\"ahler two
form on $\tilde{T}(\S )$ in two equivalent ways (see Proposition
9.1). This two form projects to the Teichm\"uller space $T(\S )$:

\vskip .2 cm

\paragraph{\bf Theorem 9.2}
{\it The Weil-Petersson two form on $\tilde{T}(\S_G )$ projects to
a non-degenerate two form on the Teichm\"uller space $T(\S_G )$.
Moreover, the two form is invariant under $Mod_{BLP}(\S )$ (see
Theorem 9.4).}

\vskip .2 cm

This paper is organized as follows.  Section~2 describes
preliminary material including profinite completion and the
universal cover of a solenoid.  Section~3 then develops the
Teichm\"uller theory of the punctured solenoid using
quasiconformal mappings, while Section~4 introduces the
representation-theoretic viewpoint of hyperbolic geometry.
Section~5 then develops the decorated theory of the solenoid and
gives the basic parametrization of decorated Teichm\"uller space
by lambda length coordinate functions.  Section 6 develops the
convex hull construction in Minkowski space and leads to a
canonical decomposition of decorated Teichm\"uller space; some of
the background material as well as other basic definitions and
results are surveyed in Appendix A, which probably should be read
before Section 6. Section~7 then introduces the mapping class
group of the punctured solenoid, presents basic results about it,
and furthermore analyzes its action on the decomposition elements
in decorated Teichm\"uller space. In Section~8, we introduce a
collection of elements of the mapping class group derived from
Whitehead moves, and by comparing their action on the decomposition
elements with that of the modular group, we conclude that,
together with $PSL_2(\mathbb{Z} )$, they generate this group.
Finally in Section~9, a two form parallel to the Weil-Petersson
K\"ahler form is introduced, and its basic properties are
described.  The short section 10 contains closing remarks and open
questions. Appendix~A surveys background material on
homeomorphisms of the circle from \cite{P1}, and Appendix~B
surveys background material on punctured surfaces of finite type
from \cite{P2} and explicates the presentation of their modular groups.

\vskip .3cm

\noindent {\bf Acknowledgements}~It is a pleasure for both authors
to thank Francis Bonahon, Bob Guralnick, Jack Milnor, Dennis
Sullivan, and Mahmoud Zeinalian for crucial comments and
questions.

\vskip .5cm

\section{Preliminaries}

\vskip .3 cm

Let $(S,x)$ be a punctured torus with a base point $x\in S-\{{\rm puncture}\}$. Consider the
set of all pointed unbranched finite coverings $\pi_i:(S_i,x_i)\to
(S,x)$ by punctured surfaces $( S_i,x_i)$ such that $\pi_i
(x_i)=x$ for all $i$. There is a natural partial ordering $\leq$ on the set of
coverings as follows. $(S_i,x_i)\leq (S_j,x_j)$ if
$\pi_j:(S_j,x_j)\to (S,x)$ factors through $\pi_i:(S_i,x_i)\to
(S,x)$, namely, if there is a pointed unbranched finite covering
$\pi_{j,i}:(S_j,x_j)\to (S_i,x_i)$ such that $\pi_j
=\pi_{j,i}\circ\pi_i$. The system of coverings $(S_i,x_i)$ is
an inverse system because given any two coverings there exists a
third covering which is greater than or equal to both.

\vskip .1 cm

\paragraph{\bf Definition 2.1} The punctured solenoid $\S$ is the
inverse limit of the system of coverings of $(S,x)$.

\vskip .1 cm

If we start the above construction from any other punctured
surface of negative Euler characteristic, we obtain a homemorphic
inverse limit. To see this, observe that any two pointed punctured
surfaces of negative Euler characteristic have a common pointed
cover. This implies that they have a common co-final subsystem of
coverings which implies that the inverse limits are homeomorphic.
The punctured solenoid $\S$ is an initial object for the
category of punctured surfaces with negative Euler characteristic with
morphisms finite unbrached covers.

$\S$ is locally homeomorphic to (2-disk)$\times$(Cantor set). The
path components are called leaves. Each leaf is homeomorphic to
the unit disk and is dense in $\S$. These observations are made
in similar fashion to the case of the inverse limit of
coverings of a compact surface of genus greater than one (Sullivan's
universal hyperbolic solenoid). However, the punctured solenoid is
not a compact space unlike that of Sullivan
\cite{Sul}. To see this, note that $\S$ is a closed subset of the
product of all covering surfaces of the punctured torus consisting
of all backward trajectories with initial points on the punctured
torus. A sequence of backward trajectories with the first
coordinates converging to the puncture on the base punctured torus
has no convergent subsequence.

Let $G$ be any finite index subgroup of $PSL_2({\mathbb
Z})$, so ${\mathbb D}/G$ is an orbifold. $G$
has characteristic subgroups
$$G_N=\cap\{ \Gamma < G:[\Gamma :G]\leq N\},$$
for each $N\geq 1$, and these are nested $G_N \leq G_{N+1}$.
Define a metric on $G$ by taking the distance between $\gamma
,\delta\in G$ to be
$${\min}\{ N^{-1}:\gamma\delta ^{-1}\in G_N\},
$$
and define the {\it profinite completion} $\hat G$ of $G$ as the
metric completion of $G$, i.e., equivalence classes of Cauchy
sequences in $G$ for the above metric.  Termwise multiplication of
Cauchy sequences gives $\hat G$ the structure of a topological
group. Moreover, $\hat{G}$ is a compact space homeomorphic to a
Cantor set (see for example \cite{Odd}). $G$ has a natural
embedding into $\hat{G}$ by mapping each $g\in G$ to
the constant Cauchy sequence $[g,g,\ldots ]\in\hat{G}$,
and the image of $G$ is dense in $\hat{G}$ by definition. Since $G$ is
naturally a subgroup of $\hat{G}$, we get a continuous right action of
$G$ on $\hat G$, i.e.,
$g\times [a_1,a_2,\ldots ]\mapsto [a_1g^{-1},a_2g^{-1},\ldots
] =[a_1,a_2,\ldots ] g^{-1}$ for $g\in G$ and
$[a_1,a_2,\dots ]\in\hat{G}$.

To alternatively define the punctured solenoid ${\mathcal H}$ as a
topological space, for definiteness fix the once-punctured torus group
$G<PSL_2({\mathbb Z})$, the unique charactersitic subgroup of index six,
let
$\hat G$ denote its profinite completion, and define the
``$G$-tagged solenoid'' $${\mathcal H}_G=({\mathbb D}\times
\hat{G})/G,$$ where $\gamma\in G$ acts on $(z,t)\in{\mathbb D}
\times\hat{G}$
by
$\gamma (z,t)=(\gamma z,t\gamma ^{-1})$. The space $\S_G$ is
homeomorphic to $\S$ following the proof of \cite{Odd} in the
compact case.

As in definition 2.1, we can start from any
finite index subgroup of $G$ and repeat the above construction to
obtain space homeomorphic to $\S_G$. The construction of $\S_G$ using a
particular group $G$ yields more structure on the leaves, namely, the
hyperbolic metric on leaves.  Different choices of Fuchsian groups give
different hyperbolic metrics on the solenoid if the groups have no common
finite index subgroups (see \cite{NS}).

The topological picture of $\S$ is easy to understand in the
model $\S_G$. Namely (in the notation of the introduction), let $\gamma
,\delta\in G$ be the hyperbolic generators of $G$ which map
$e_{1/2}$ and $e_{2/1}$ onto $e_{-2/1}$ and $e_{-1/2}$,
respectively. Let $Q$ be the quadrilateral in ${\mathbb D}$ with
boundary sides $e_{1/2},e_{2/1},e_{-2/1}$ and $e_{-1/2}$. Thus, $Q$
is a fundamental polygon for $G$, and $g$ and $h$ are side
identifications of $Q$. The set $Q\times\hat{G}$ is a fundamental
set for the action of $G$ on ${\mathbb D}\times\hat{G}$. The side
$e_{1/2}\times t$ is identified with $e_{-2/1}\times t\gamma^{-1}$
and the side $e_{2/1}\times t$ is identified with $e_{-1/2}\times
t\delta^{-1}$, for all $t\in\hat{G}$. These are the only
identifications by $G$ on the fundamental set $Q\times\hat{G}$.
Therefore, the punctured solenoid is obtained by ``sewing'' the
Cantor set of polygons along their boundary sides according to the
action of $G$ on $\hat{G}$. This immediately shows that leaves are
unit disks and that each leaf is dense (by the density of $G$ in
$\hat{G}$).

\section{Quasiconformal definition of the Teichm\"uller space $T
(\S_G )$}

\vskip .5cm

In this section, we define the Teichm\"uller space of the punctured
solenoid
$\S_G \equiv {\mathbb D}\times_G\hat{G}$ in the spirit of Ahlfors-Bers
theory \cite{AB} and  in such fashion that the union of the natural lifts
of Teichm\"uller spaces of all finite sheeted covers of the modular curve
is dense. The condition of continuity for the transverse variation, which
is familiar from the compact case
\cite{Sul} and \cite{Sa}, requires an extra stipulation in its
definition in the punctured case due to the non-compactness.

Fix a metric (which is complete with constant curvature -1) on
$\S_G$ coming from the lift of the hyperbolic metric on the
modular curve (or the covering of the modular curve by the
punctured torus with the covering group $G<PSL_2\mathbb{Z}$). This
is a transversely locally constant metric, namely there exists a
subatlas of $\S_G$ with the metric being constant for the
transverse variation in the charts of the subatlas. Let $\S_d$ be
the differentiable structure on the solenoid compatible with the
hyperbolic metric on $\S_G$.

\vskip .3cm

\paragraph{\bf Definition 3.1}
A {\it hyperbolic metric} on $\S_d$ is an assignment of a hyperbolic
metric (a complete metric of curvature $-1$) on each leaf that
varies continuously (in the $C^{\infty}$-topology) for the
transversal variations in the local charts.

\vskip .3 cm

A hyperbolic metric on $\S_d$ gives a complex structure on $\S_d$.
The complex structure varies continuously for the transverse
variations on the local charts. The converse also holds \cite{Ca}.

Let $\S$ be an arbitrary complex solenoid with the underlying
differentiable solenoid $\S_d$. As a preliminary for the
definition of the Teichm\"uller space, we define a
(differentiable) quasiconformal map from the fixed solenoid $\S_G$
onto $\S$. Recall that a leaf of the solenoid intersects any local
chart infinitely many times. Each component of the intersection is
called a {\it local leaf} of the solenoid.

A chart for $\S_G$ is modeled on an open set in (hyperbolic
disk)$\times$(transverse Cantor set). Since the natural projection
$\pi :{\mathbb D}\times\hat{G}\to ({\mathbb D}\times\hat{G})/G$ is a local
homeomorphism, there exists $U\subset{\mathbb D}$ such that the set
$U\times\hat{G}$ together with the map $\pi
:U\times\hat{G}\to\S_G$ is a local chart for $\S_G$. In
particular, we choose $U\subset{\mathbb D}$ to be a hyperbolic disk
with center $0\in{\mathbb D}$ and radius  smaller than the
injectivity radius of ${\mathbb D}/G$. The chart map $\pi
:U\times{\mathbb D}\to\S_G$ is an isometry on the leaves, and the
natural identification $U\times t\equiv U\times t_1$ is also an
isometry, for all $t,t_1\in\hat{G}$. For two fixed global leaves of
$\S_G$ and a choice of their corresponding local leaves in the chart
$U\times\hat{G}$, there is thus a preferred isometric identification of
the two local leaves. This gives an isometric identification
between two open sets of the global leaves, and this identification
extends uniquely to a preferred isometric identification of the two global leaves
themselves.

\vskip .3cm

\paragraph{\bf Definition 3.2}
A homeomorphism $f:\S_G\to\S$ is said to be a {\it
(differentiable) }$K${\it -quasiconformal map} if its restriction
to each leaf is a $K$-quasiconformal $C^{\infty}$-map, if the
restriction of $f$ to any local chart varies continuously for the
transversal variation in the $C^{\infty}$-topology, and if $f$
varies continuously (in the quasiconformal topology) for the
transversal variation on global leaves (with two leaves identified
using their local leaves in the fixed chart as above).

\vskip .3 cm

\paragraph{\bf Remark 3.3} The quasiconformal topology is induced by
the semi-norm on the space of quasiconformal maps: the distance
between two quasiconformal maps $f$ and $g$ is induced by the sup norm
of the Beltrami differential of the map $f\circ g^{-1}$.
(This is
only a semi-norm because the Beltrami coefficient of a conformal map
vanishes identically.) The quasiconformal topology induces a topology on
the quotient of the space of quasiconformal maps by post-composition with
M\"obius transformations, and the resulting space is Hausdorff.

\vskip .3 cm

\paragraph{\bf Remark 3.4} The natural requirement of continuity in the
$C^{\infty}$-topology for the transversal variations in the local
charts is given by Sullivan \cite{Sul} (see also \cite{Sa}) for
compact Riemann surface laminations although other transverse
structures are also possibly interesting.  In the case
of the universal hyperbolic compact solenoid (see \cite{Sul},
\cite{NS} and \cite{Sa}), this is enough to guarantee that the
union of Teichm\"uller spaces of compact surfaces is dense in the
Teichm\"uller space of the compact solenoid. However, for the
punctured solenoid we require the additional quasiconformal continuity because of the
non-compactness.

\vskip .3cm

To make sense of the identifications of the
different global leaves, we are forced to start from a locally
constant hyperbolic structure to make a proper choice of the
identifications of the global leaves.

\vskip .3cm

\paragraph{\bf Definition 3.5}
The {\it Teichm\"uller space} $T (\S_G)$ of the punctured solenoid
$\S_G$ is the space of all differentiable quasiconformal maps
$f:\S_G\to \S$, where $\S$ is an arbitrary hyperbolic solenoid
with the underlining differentiable solenoid $\S_d$ up to
the following
equivalence relation: two maps $f_1:\S_G\to\S_1$
and $f_2:\S_G\to\S_2$ are {\it Teichm\"uller equivalent} if there is
a hyperbolic isometry $c:\S_1\to \S_2$ such that $f_2^{-1}\circ c\circ
f_1:\S_G\to\S_G$ is homotopic to the identity by a bounded
homotopy with respect to the hyperbolic metric.

\vskip .3cm

We already introduced the universal covering of $\S_G$ by the natural
projection map $\pi :{\mathbb D}\times\hat{G}\to ({\mathbb
D}\times\hat{G})/G$. Note that the restriction of $\pi$ to each leaf of
the covering is an isometry for the hyperbolic metrics of the leaf in the
covering and the leaf of $\S_G$. However, countably many leaves of the
covering are mapped onto a single leaf of $\S_G$. Given an arbitrary
hyperbolic punctured solenoid
$\S$, there exists a quasiconformal map $f:\S_G\to\S$. We
next introduce a universal covering of $\S$ using the existence of $f$,
similarly to \cite{Sa}.

Let $U\times\hat{G}$, $0\in U$, be a local chart of $\S_G$ such
that $f|_{(\pi (U\times\hat{G} ))}$ is a homeomorphism. Choose a
local chart $V\times T\to\S$ for $\S$ which contains $f(\pi
(0\times\hat{G}))$ and let $\psi :V\times T\to\S$ be isometric
local chart map. Identify $T$ with $\hat{G}$ via the
correspondence $\hat{G}\to 0\times\hat{G}\to f(0\times\hat{G})\to
T$. Finally, define the universal covering $\pi_{\S} :{\mathbb
D}\times\hat{G}\to\S$ as follows. Choose a continuous isometric embedding
$i:V\times \hat{G}\to {\mathbb D}\times \hat{G}$ such that
$i\circ\psi^{-1} (f(0\times\hat{G}))=0\times\hat{G}$ and define
$\pi_{\S}$ by extending the local isometry $i\circ\psi^{-1}$ from
$V\times \hat{G}\subset {\mathbb D}\times \hat{G}$ to the global
leafwise isometry from ${\mathbb D}\times \hat{G}$ to $\S_{G}$. Note
that $f:\S_G\to\S$ lifts to a quasiconformal map $\tilde{f}:{\mathbb
D}\times\hat{G}\to {\mathbb D}\times\hat{G}$ by the formula $\tilde{f}
(z,t)=\pi_{\S}^{-1}( f\circ\pi (z,t),t)$, for $z\in{\mathbb D}$ and
$t\in\hat{G}$. Further, the action of any $\gamma\in G$ on ${\mathbb
D}\times\hat{G}$ is conjugated by $\tilde{f}$ to the action on the
universal cover of $\S$. This conjugated action is also an isometry from
one leaf onto another, but it depends on the leaf, i.e., it is not
constant in $t\in \hat{G}$ (this follows from the formula for $\tilde{f}$,
or see
\cite{Sa} for a similar argument). In the proposition below, we construct
the hyperbolic punctured solenoid by constructing first the universal
covering and the covering group from a Beltrami coefficient, thus
reversing the above construction.

The leafwise Beltrami coefficient $\mu$ on $\S_G$ of a
differentiable quasiconformal map $f:\S_G\to \S$ is smooth,
$\|\mu\|_{\infty}\leq k<1$, it varies continuously for the
transversal variations in the smooth topology and
$\|\mu_{t_1}-\mu_t\|_{\infty}\to 0$ as $t_1\to t$, where
$\mu_{t_1},\mu_t$ are the restrictions of $\mu$ to the global
leaves $t_1,t$ transversely identified as above with the unit disk
${\mathbb D}$. Conversely, we show

\vskip .3cm

\paragraph{\bf Proposition 3.6} {\it Let $\mu$ be a Beltrami
coefficient on $\S_G$ so that $\mu$ is smooth, has $L^\infty$-norm
bounded above by $k<1$, varies continuously for transverse
variations in the smooth topology, and so that
$\|\mu_{t_1}-\mu_t\|_{\infty}\to 0$ as $t_1\to t$ in the
transverse direction. Then there exists a hyperbolic punctured
solenoid $\S$ and a (differentiable) quasiconformal map
$f:\S_G\to\S$ such that the Beltrami coefficient of $f$ is $\mu$.
Moreover, $f$ and $\S$ are unique up to post-composition with a
hyperbolic isometry.}

\vskip .3cm

\paragraph{\bf Proof}
The natural quotient map $\pi :{\mathbb D} \times\hat{G}\to {\mathbb D}
\times_G\hat{G}= \S_G$ is the universal covering map of
$\S_G$. Denote by $\tmu$ the lift of $\mu$ to the universal
covering; hence $\tmu$ satisfies

\begin{equation}
\label{invbelt} \tmu (\gamma
z,t\gamma^{-1})~{\overline{\gamma^{'}(z)}}=\tmu
(z,t)~{\gamma^{'}(z)}
\end{equation}
for all $\gamma\in G$.

Thus, $\tmu (\cdot ,t)$, for $t\in\hat{G}$, varies continuously
with $t$ in the norm $\|\cdot\|_{\infty}$ on ${\mathbb D}\equiv{\mathbb
D}\times\{ t\}$ by definition of $\mu$. We solve the Beltrami equation for
$\tmu (\cdot ,t)$ on each ${\mathbb D} \times\{ t\}$ such that the
solution
$\tilde{f}_t:{\mathbb D} \times\{ t\}\to {\mathbb D} \times\{ t\}$ fixes
$1$, $-1$ and $-i$.

The resulting $\tilde{f}_t$ is smooth and varies continuously in the
smooth and quasiconformal topology by the corresponding properties of
$\tmu$ (see \cite{AB}). By (\ref{invbelt}), the maps $\tilde{f}_t$
conjugate the action of $\gamma\in G$ on ${\mathbb D} \times\hat{G}$ into
the action of $\gamma_{\mu}\in G_{\mu}$ on ${\mathbb D} \times\hat{G}$,
where each $\gamma_{\mu}$ acts by isometry leafwise and varies
continuously for the transverse variation. Thus, $\tilde{f}$
projects to a homeomorphism $f:{\mathbb D} \times\hat{G}/G\to {\mathbb D}
\times\hat{G}/G_{\mu}$ which is a differentiable quasiconformal map
from $\S_G$ to the hyperbolic punctured solenoid $\S = {\mathbb D}
\times\hat{G}/G_{\mu}$. (For a similar construction, see
\cite[Section 2]{Sa}.) $\Box$

\vskip .3cm

We thus obtain an equivalent definition of $T(\S_G)$:

\vskip .3cm

\paragraph{\bf Definition 3.7} The {\it Teichm\"uller space} $T(\S_G )$
is the space of all Beltrami coefficients on $\S_G$ which satisfy
the conditions of Proposition 3.1 up to an equivalence. Two
Beltrami coefficients are {\it equivalent} if their corresponding
quasiconformal maps are equivalent.

\vskip .3cm

Sullivan \cite{Sul} showed that union of the Teichm\"uller spaces of all
compact surfaces is dense in the Teichm\"uller space of the
compact solenoid, and we next prove the corresponding result for the
punctured solenoid. Any punctured surface of finite area is covered by the
punctured solenoid, and in fact, this covering is a principal fiber
bundle. The Teichm\"uller space of a punctured surface lifts to $T (\S_G
)$ by lifting hyperbolic metrics on the surface to locally constant
hyperbolic metrics on $\S_d$. For example, the hyperbolic metric
on $\S_G$ is the lift of the hyperbolic metric on the modular
curve ${\mathbb D}/PSL_2({\mathbb Z})$, which in turn lifts to a
hyperbolic metric on the punctured torus or any other surface covering
the modular curve.

\vskip .3cm

\paragraph{\bf Theorem 3.8}
{\it The union of natural lifts of Teichm\"uller spaces of all
punctured hyperbolic surfaces is dense in the Teichm\"uller space
$T(\S_G )$.}

\vskip .3cm

\paragraph{\bf Proof}
It is enough to show that for any Beltrami coefficient $\mu$ on
$\S_G$ there is a sequence of locally constant Beltrami
coefficients $\mu_n$ which approximate $\mu$ in the Teichm\"uller
topology.

Let $\tmu$ be the lift of $\mu$ to ${\mathbb D}\times\hat{G}$, so $\tmu$
satisfies (\ref{invbelt}). Let $G_n$ be the intersection of all
index at most $n$ subgroups of G (see Section 2). Let $P$ be a
fundamental polygon for the action of $G_n$ on ${\mathbb D} $. We define
$\tmu_n(z,t)=\tmu (z,id)$ for $z\in P$ and $t\in \hat{G}_n$, and
extend $\tmu_n$ to ${\mathbb D} \times\hat{G}_n$ by the action of $G_n$.
Since ${\mathbb D} \times\hat{G}/G\equiv {\mathbb D} \times\hat{G}_n/G_n$
and
$\tmu_n$ is close to $\tmu$ in the Teichm\"uller topology on each
leaf, we obtain the required locally constant sequence
approximating $\tmu$. We may replace $\tmu_n$, for instance with
the barycentric extension of its boundary values (see Douady-Earle
\cite{DE}), to produce a Teichm\"uller equivalent Beltrami
coefficient which is smooth and transversely continuous for both
the smooth and quasiconformal topologies. $\Box$

\section{The representation-theoretic definition of $T(\S_G )$}

We give an alternative definition of the Teichm\"uller space using
the universal covering and the representation of the covering
group of the solenoid. This definition is motivated by the
finite-dimensional Teichm\"uller theory of punctured
surfaces. We use the universal covering construction above which
was adopted from \cite{Sa} to the case of the punctured solenoid.

\vskip .3cm

Let us consider the collection $Hom(G\times\hat{G},PSL_2\RR )$ of
all functions $\rho :G\times\hat{G}\to PSL_2\RR $ satisfying the
following three properties:

\vskip .3cm

\leftskip .3cm

\noindent {Property 1:} $\rho$ is continuous;

\vskip .3cm

\noindent {Property 2 [$G$-equivariance]:} for each $\gamma _1,\gamma
_2\in G$ and
$t\in \hat {G}$, we have
$$
\rho(\gamma _1\circ\gamma
_2,t)=\rho(\gamma _1,t\gamma _2^{-1})\circ \rho(\gamma _2,t);
$$

\vskip .3cm

\noindent {Property 3:} for every $t\in\hat{G}$, there is a
quasiconformal mapping $\phi _t:{{\mathbb D}}\to{{\mathbb D}}$ depending
continuously on $t\in\hat{G}$ so that for every $\gamma\in G$, the
following diagram commutes:

\vskip .3cm

\begin{equation*}
\begin{array}{lll}
(z,t)\ \ \ \ \ \ \ \ \ \ \ \ \  & \mapsto\ \ \ \ \ \ \ \ \ \
 & (\gamma z,t\gamma ^{-1})\\
\\
{{\mathbb D}}\times\hat{G} & \to & {{\mathbb D}}\times\hat{G}\\
\\
\phi_t\times{\rm id}\downarrow & & \downarrow \phi_{t\gamma
^{-1}}\times{\rm id} \\
\\
{{\mathbb D}}\times\hat{G} & \to & {{\mathbb D}}\times\hat{G} \\
\\
(\phi_t(z),t) & \mapsto & (\rho(\gamma
,t)\circ\phi_t(z)=\phi_{t\gamma^{-1}}\circ\gamma (z),t\gamma
^{-1})
\end{array}
\end{equation*}

\leftskip=0ex

\vskip .3cm

\noindent  As to property 1, notice that since $G$ is discrete,
$\rho$ is continuous if and only if it is so in its second
variable only. Property 2 is a kind of homomorphism property of
$\rho$ mixing the leaves; notice in particular that taking $\gamma
_2=~I$ gives $\rho(I,t)=I$ for all $t\in\hat{G}$.  Property 3
mandates that for each $t\in \hat{G}$, $\phi _t$ conjugates the
standard action of $\gamma\in G$ on ${{\mathbb D}}\times\hat{G}$ at the
top of the diagram to the action
$$\gamma _\rho:(z,t)\mapsto (\rho (\gamma ,t)z,t\gamma ^{-1})$$ at
the bottom, and we let $G_\rho=\{ \gamma _\rho:\gamma\in
G\}\approx G$. Notice that the action of $G_\rho$ on ${{\mathbb
D}}\times\hat{G}$ extends continuously to an action on
$S^1\times\hat{G}$.  We finally define the solenoid (with marked
hyperbolic structure)
$$
\S_\rho =({{\mathbb D}}\times _\rho\hat{G})=({{\mathbb D}}\times
\hat{G})/G_\rho.
$$

\vskip .3cm

\noindent Define the group $Cont(\hat{G},PSL_2\RR )$ to be the
collection of all continuous maps $\alpha:\hat{G}\to PSL_2\RR $,
where the product of two $\alpha ,\beta\in Cont(\hat{G},PSL_2\RR
)$ is taken pointwise $(\alpha\beta) (t)=\alpha (t)\circ\beta (t)$
in $PSL_2\RR $. $\alpha\in Cont(\hat{G},PSL_2\RR )$ acts
continuously on $\rho\in Hom(G\times \hat{G},PSL_2\RR )$ according
to
$$(\alpha\rho)(\gamma, t)=\alpha(t\gamma^{-1})
\circ\rho(\gamma ,t)\circ \alpha^{-1} (t).$$

\vskip .3cm

We introduce the topology on $Hom(G\times\hat{G},PSL_2\RR )$.
Consider the natural metric $d$ on $PSL_2\RR$ induced by
identifying it with the unit tangent bundle of the unit disk ${\mathbb
D}$. Let $\rho_1,\rho_2\in Hom(G\times\hat{G},PSL_2\RR )$ and let
$\gamma_1,\ldots ,\gamma_j\in G$ be a generating set of $G$.
The distance between $\rho_1$ and $\rho_2$ is given by
\begin{equation}
\label{repmet} \max_{1\leq i\leq j,\
t\in\hat{G}}d(\rho_1(\gamma_i,t), \rho_2(\gamma_i,t)).
\end{equation}
This metric is not canonical, but any such two metrics induce the
same topology.

The topology on $Hom'(G\times\hat{G},PSL_2\RR
)=Hom(G\times\hat{G},PSL_2\RR )/Cont(\hat{G},PSL_2\RR )$ is the
quotient topology of the above topology on
$Hom(G\times\hat{G},PSL_2\RR )$.

\vskip .3cm

\paragraph{\bf Theorem 4.1}{\it There is a natural homeomorphism of
the Teichm\"uller space of the solenoid ${\S_G}$ with
$$Hom'(G\times\hat{G},PSL_2\RR
)=Hom(G\times\hat{G},PSL_2\RR )/Cont(\hat{G},PSL_2\RR ).$$}

\vskip .3cm

\paragraph{\bf Proof} Let $\mu$ be a Beltrami coefficient on
$\S_G$ which represents a point in $T(\S_G )$. Let $\tmu$ be the
lift of $\mu$ to ${\mathbb D}\times \hat{G}$ and let
$\tilde{f}:{\mathbb D}\times\hat{G}\to{\mathbb D}\times\hat{G}$ be the
leafwise solution of the corresponding Beltrami equation which fixes $1$,
$-1$ and $-i$ on the boundary of each leaf, as in Proposition 3.6.
Since $\tmu$ satisfies (\ref{invbelt}), then for each $\gamma\in
G$ there exists $\gamma_{\mu}:{\mathbb D}\times \hat{G}\to{\mathbb
D}\times
\hat{G}$ which is leafwise a M\"obius transformation such that
\begin{equation}
\label{three}
(\tilde{f}\circ\gamma
)(z,t)=(\gamma_{\mu}\circ\tilde{f})(z,t)
\end{equation}
for all $(z,t)\in {\mathbb D}\times \hat{G}$ and $\gamma\in G$. We note
that $\tilde{f}$ fixes each leaf of ${\mathbb D}\times \hat{G}$, and that
$\gamma$ and $\gamma_{\mu}$ exchange leaves in the same
fashion ${\mathbb D}\times\{ t\}\mapsto{\mathbb D}\times\{
t\gamma^{-1}\}$.

We define the representation $\rho_{\mu}:G\times\hat{G}\to
PSL_2\RR$ by
\begin{equation}
\label{repres}
\rho_{\mu}(\gamma ,t)=\gamma_{\mu}(\cdot ,t)
\end{equation}
for $\gamma\in G$ and $t\in\hat{G}$.

We check that properties 1, 2 and 3 hold. To see that $\rho_{\mu}$
is continuous, note that an element of $PSL_2\RR$ is uniquely
determined by specifying the image of three arbitrary points on
$S^1$, and that it depends continuously on this image. Therefore,
Property 1 holds because $\tilde{f}$ varies continuously on
$S^1\times\hat{G}$ in the $C^0$-topology.

Property 2 follows from the computation:
\begin{eqnarray*}
&\rho_{\mu}(\gamma_1\circ\gamma_2
,t)=\tilde{f}\circ\gamma_1\circ\gamma_2\circ\tilde{f}^{-1} (\cdot
,t)= \tilde{f}\circ\gamma_1\circ\tilde{f}^{-1}\circ
\tilde{f}\circ\gamma_2\circ\tilde{f}^{-1} (\cdot ,t)=\\
& [(\tilde{f}\circ\gamma_1\circ\tilde{f}^{-1})(\cdot
,t\gamma_2^{-1})]\circ
[(\tilde{f}\circ\gamma_2\circ\tilde{f}^{-1}) (\cdot
,t)]=\rho_{\mu}(\gamma_1 ,t\gamma_2^{-1})\circ \rho_{\mu}(\gamma_2
,t).
\end{eqnarray*}
Thus, $\phi_t:=\tilde{f}(\cdot ,t)$ is a quasiconformal map
which satisfies Property 3.

Upon taking quotients, the assignment (\ref{repres}) induces a
well-defined map from $T(\S_G )$ to $Hom'(G\times\hat{G},PSL_2\RR
)=Hom(G\times\hat{G},PSL_2\RR )/Cont(\hat{G},PSL_2\RR )$. To see
this, simply note that $f_t|_{S^1}$ up to post-composition with an
element of $PSL_2\mathbb{R}$ depends only on the Teichm\"uller
class of the Beltrami coefficient $\mu$. Thus, $\gamma_{\mu}$
is determined up to the same ambiguity, and the map is well-defined
on the quotients.

We prove the continuity of this map at an arbitrary $[\mu_1]\in
T(\S_G )$, where $[\mu_1 ]$ denotes the Teichm\"uller class of the
Beltrami coefficient $\mu_1$. Let $\mu_2$ be a Beltrami
coefficients representing a points in $T(\S_G )$ in a small
neighborhood of $[\mu_1]$ such that $\|\mu_1-\mu_2\|_{\infty}\to
0$ as $[\mu_2]\to [\mu_1]$. Their corresponding maps $\tilde{f}_1$
and $\tilde{f}_2$ are thus close as quasiconformal maps for each
$t\in\hat{G}$, and they both fix $1$, $-1$ and $i$ on each leaf.

Recall that the group $G<PSL_2\mathbb{Z}$ for the punctured torus
is generated by two hyperbolic elements $\gamma_1$ and $\gamma_2$
which carry $e_{1/2}$ onto $e_{-2/1}$ and $e_{2/1}$ onto
$e_{-1/2}$, respectively.

It follows that $G_{\mu_k}$ are generated by their
conjugates
$\gamma_1^{\mu_k},\gamma_2^{\mu_k}$, for $k=1,2$. To estimate
$d(\gamma_1^{\mu_1}(\cdot ,t),\gamma_1^{\mu_2}(\cdot ,t))$ and
$d(\gamma_2^{\mu_1}(\cdot ,t),\gamma_2^{\mu_2}(\cdot ,t))$, it is
enough to compare the corresponding images of $-1$, $-i$ and $1$,
for all $t\in\hat{G}$. To this end, we use the formula
$$
\gamma_1^{\mu_k}(\cdot ,t)=\tilde{f}_k(\cdot ,t\gamma_1^{-1})\circ
\gamma_1(\cdot ,t)\circ \tilde{f}_k(\cdot , t)^{-1},
$$
for all $t\in\hat{G}$, and similarly for $\gamma_2^{\mu_k}$. Thus,
$\gamma_1^{\mu_k}(-1 ,t)=\tilde{f}_k(i,t\gamma_1^{-1})$,
$\gamma_1^{\mu_k}(-i ,t)=\tilde{f}_k(1,t\gamma_1^{-1})$ and $
\gamma_1^{\mu_k}(1 ,t)=\tilde{f}_k(\gamma_1(1),t\gamma_1^{-1}). $
Therefore, it is enough to show that pairs
$\tilde{f}_1(i,t\gamma_1^{-1})$ and
$\tilde{f}_2(i,t\gamma_1^{-1})$, $\tilde{f}_1(1,t\gamma_1^{-1})$
and $\tilde{f}_2(1,t\gamma_1^{-1})$, and
$\tilde{f}_1(\gamma_1(1),t\gamma_1^{-1})$ and
$\tilde{f}_2(\gamma_1(1),t\gamma_1^{-1})$ are close in the angle
metric on $S^1$. This follows since $\tilde{f}_1$ and
$\tilde{f}_2$ are properly normalized and their Beltrami
coefficients are close in the supremum norm uniformly in
$t\in\hat{G}$. The similar statement holds for $\gamma_2^{\mu_k}$.
The continuity is proved.

Assume that $\rho_{\mu_1}=\rho_{\mu_2}=\rho$ for two Beltrami
coefficients $\mu_1,\mu_2$ representing elements of $T(\S_G )$. We
need to show that $\tilde{f}_1(\cdot ,t)|_{S^1}=\tilde{f}_2(\cdot
,t)|_{S^1}$ for all $t\in\hat{G}$. This implies that there is a
bounded homotopy through quasiconformal maps between $\tilde{f}_1$
and $\tilde{f}_2$ which is invariant under $G_{\rho}$ (the proof
for the compact solenoid \cite{MS} extends directly to the
punctured solenoid), which says that $\mu_1$ and $\mu_2$ are
Teichm\"uller equivalent, so the map is one to one.

It is enough to show that $\tilde{f}_1(x,t)=\tilde{f}_2(x,t)$ for
all $x\in\bar{\mathbb{Q}}$ and $t\in\hat{G}$. Note that the
equality holds at $1/0,0/1,1/1\in\bar{\mathbb{Q}}$ by our
normalization. Furthermore, at least one edge of the triangle with  vertices
$1/0,0/1,1/1$ is mapped onto any other edge of the Farey
tesselation by an element of the once-punctured torus group $G$, where $G$ is generated by two
hyperbolic elements $\gamma_1$ and $\gamma_2$. We proceed
inductively. Namely, assume that $x\in\bar{\mathbb{Q}}$ is one
vertex of edge $e\in\tau_{*}$ and that at vertices of
$e_1\in\tau_{*}$ the equality of the two maps hold for all
$t\in\hat{G}$ and that $\gamma_i (e_1)=e$ for either $i=1$ or
$i=2$. Let $x_1$ be the vertex of $e_1$ such that
$\gamma_i(x_1)=x$. By (\ref{three}), it follows that
$\tilde{f}_j(x,t)=(\gamma_i^{\mu_j}\circ\tilde{f}_j)(x_1,t\gamma_i
)$ for $j=1,2$. By our assumption, $\gamma_i^{\mu_1}\equiv
\gamma_i^{\mu_2}$ and $\tilde{f}_1(x_1,t\gamma_i
)=\tilde{f}_2(x_1,t\gamma_i )$ for all $t\in\hat{G}$. Thus,
$\tilde{f}_1(x,t)=\tilde{f}_2(x,t)$ for all $t\in\hat{G}$. The
inductive step is complete, and the injectivity follows.

The map is surjective because ${\mathbb D}\times\hat{G}/G_{\rho}$ is a
quasiconformal image of $\S_G$. The quasiconformal map is the
projection of the family $\phi_t:{\mathbb D}\to{\mathbb D}$, for
$t\in\hat{G}$, to the quotients $\S_G$ and $\S_{\rho}$ obtained by the
actions of
$G$ and $G_{\rho}$ on ${\mathbb D}\times\hat{G}$.

We must finally show that the inverse map is continuous. Suppose that
$\rho_1$ and $\rho_2$ are two representations of the group $G$
which are close in the metric (\ref{repmet}). Let $\phi_1,\phi_2
:{\mathbb D}\times\hat{G}\to
{\mathbb D}\times\hat{G}$ be conjugating maps for $\rho_1,\rho_2$
respectively, from Property 3. It is enough to find another pair
of conjugating maps $\tilde{f}_1,\tilde{f}_2$ such that the
quasiconformal constant of $\tilde{f}_2\circ\tilde{f}_1^{-1}$
converges to $1$ as $\rho_1$ converges to $\rho_2$.

Let $Q$ be the ideal hyperbolic rectangle inside ${\mathbb D}$ with
vertices $1$, $i$, $-1$ and $-i$, so $Q\times\hat{G}$ is a
fundamental domain for the action of the punctured torus group $G$
on ${\mathbb D}\times\hat G$. Let $\gamma_1,\gamma_2\in G$ be fixed
hyperbolic generators of the punctured torus group $G$ as above.
The sides of $Q\times t$ are identified to the sides of $Q\times
t\gamma_1^{-1}$ and $Q\times t\gamma_2^{-1}$, for $t\in\hat{G}$.
Note that ${\mathbb D}\times \hat{G}/G\equiv Q\times\hat{G}/G\equiv \S_G$.

Let $G_{k}=\rho_k(G)$ and let $Q_k^t$ be the ideal rectangle with
vertices $-1$, $-i$, $1$ and $i_k:=\phi_k(i,t)$, for $k=1,2$ and
$t\in\hat{G}$. Thus, $F_k=\cup_{t\in\hat{G}}Q_k^t\times \{ t\}$ is
a fundamental set for $G_{k}$.

It is enough to find a quasiconformal map $\tilde{f}$ between
$F_1$ and $F_2$ with small quasiconformal constant such that
$\tilde{f}=I$ on the geodesics $e_{1/2}\times t$ and
$e_{2/1}\times t$,
$\tilde{f}=\rho_2(\gamma_1,t)\circ\rho_1(\gamma_1,t)^{-1}$ on
$\phi_1 (e_{-2/1}, t\gamma_1^{-1})$ and
$\tilde{f}=\rho_2(\gamma_2,t)\circ\rho_1(\gamma_2,t)^{-1}$ on
$\phi_1 (e_{-1/2},t\gamma_2^{-1})$. By assumption, the covering
maps $\rho_1(\gamma_j,t)$ and $\rho_2(\gamma_j,t)$, $j=1,2$ are
close and the rectangles $Q_1^t$ and $Q_2^t$ are almost equal, so
there exists a quasiconformal map $\tilde{f}$ with small
dilatation which satisfies above boundary conditions. Furthermore,
$\tilde{f}:F_1\to F_2$ lifts and extends to a quasiconformal
self-map of ${\mathbb D}\times\hat{G}$ using the actions of the covering
groups $G_{k}$, for $k=1,2$. Finally, taking $\tilde{f}_1=\phi_1$
and $\tilde{f}_2=\tilde{f}\circ\phi_1$, the continuity of the
inverse map follows. $\Box$

\section{The decorated Teichm\"uller space $\tilde{T}(\S_G)$}

Let $\S_{\rho}$ be a hyperbolic solenoid obtained from the
representation $\rho$ with corresponding quasiconformal map $\phi
:\S_G\to\S_{\rho}$. Fix a leaf of $\S_{\rho}$ and fix a point $p$
on that leaf. Consider a geodesic ray in the hyperbolic metric of
the leaf starting at $p$. If a ray leaves every compact subset of
$\S_{\rho}$ then it determines a ``puncture'' of the solenoid $\S_\rho$.
More precisely, a {\it puncture} of a hyperbolic solenoid $\S_{\rho}$ is
an equivalence class of wandering rays from points of the leaf, where two
rays are equivalent if they are asymptotic.

We may describe the punctures of $\S_{\rho}$ using the
representation $\rho$. The quasiconformal map $\phi
:{\mathbb D}\times\hat{G}\to {\mathbb D}\times\hat{G}$ extends
continuously to a leaf-wise quasi-symmetric map $\phi :S^1\times\hat{G}\to
S^1\times\hat{G}$. Recall that $\bar{\mathbb Q}\subset S^1$
parametrizes the endpoints of the standard triangulation of ${\mathbb D}$
invariant under $PSL_2\RR$. We say that a point $(p,t)\in
S^1\times\hat{G}$ is a {\it $\rho$-puncture} if $\phi
^{-1}(p,t)\in \bar{\mathbb Q}$, and a {\it puncture} of $\S_\rho$
itself is a $G_\rho$-orbit of $\rho$-punctures.  A {\it
$\rho$-horocycle} at a $\rho$-puncture $(p,t)$ is the horocycle in
${{\mathbb D}}\times\{ t\}$ centered at $(p,t)$ and a {\it
horocycle} on $\S_{\rho}$ is a $G_{\rho}$-orbit of
$\rho$-horocycles.

\vskip .3cm

 \noindent{\bf Definition 5.1}~A {\it decoration} on $\S_\rho$, or a ``decorated
hyperbolic structure'' on $\S_{\rho}$, is a function
$\tilde\rho:G\times\hat{G}\times\bar{\mathbb Q}\to PSL_2\RR \times
L^+$, where
$$\tilde\rho(\gamma ,t,q)=\rho(\gamma ,t)\times h(t,q)$$ with
$\rho(\gamma ,t)\in Hom(G\times\hat{G},PSL_2\RR )$, which
satisfies the following conditions:

\vskip .2in

\leftskip .3in

\noindent {Property 4:}~for each $t\in\hat{G}$, the image
$h(t,\bar{\mathbb Q})\subseteq L^+$ is discrete and the center of the
horocycle $h(t,q)$ is $\phi_t(q)$, for all
$(t,q)\in\hat{G}\times\bar{Q}$ (using here the identification of $L^+$ with the space of
horocycles as in Appendix A);

\vskip .1in

\noindent {Property 5:}~for each $q\in\bar{\mathbb Q}$, the
restriction $h(\cdot ,q):\hat{G}\to L^+$ is a continuous function
from $\hat{G}$ to $L^{+}$;

\vskip .1in

\noindent {Property 6:}~$h(t,q)$ is $\rho$ invariant in the sense
that $$ \rho(\gamma ,t)(h(t,q))=h(t\gamma^{-1},\rho (\gamma
,t)q).$$

\vskip .2in

\leftskip=0ex

\noindent {\bf Remark} Property 4 implies that $h(t,\bar{\mathbb Q})$ is
radially dense in $L^{+}$ because $\phi_t(\bar{Q})$ is dense in
$S^1$. The continuity from property 5 and the invariance under the
covering group from property 6 imply continuity of the map $h(\cdot
,\bar{Q})$ from $\hat{G}$ to the space of discrete and radially
dense subset of $L^{+}$, where the topology is the Hausdorff
topology on closed subsets of $L^{+}$.

\vskip .2in

Let $Hom(G\times\hat{G}\times\bar{\mathbb Q},PSL_2\RR\times L^+)$
denote the space of all decorated hyperbolic structures satisfying
the properties above. We define a topology on
$Hom(G\times\hat{G}\times\bar{\mathbb Q},PSL_2\RR\times L^+)$. A
neighborhood of $\tilde\rho(\gamma ,t,q)=\rho(\gamma ,t)\times
h(t,q)$ consists of all $\tilde\rho_1(\gamma ,t,q)=\rho_1(\gamma
,t)\times h_1(t,q)$ such that $\rho_1$ belongs to a chosen
neighborhood of $\rho$ in $Hom(G\times\hat{G},PSL_2\RR )$, and the
maps $h_1(\cdot ,q ):\hat{G}\to L^{+}$ and $h(\cdot ,q):\hat{G}\to
L^{+}$ are close in the supremum norm, for each $q\in\bar{Q}$. The
above condition and the invariance Property 6 implies that the set
$h_1(t, \bar{Q})$ is close to the set $h(t,\bar{Q})$ in the
Hausdorff metric, for each $t\in\hat{G}$.

\vskip .2in

We define the {\it decorated Teichm\"uller space} as the quotient
$$\tilde T(\S_G )=Hom(G\times\hat{G}\times\bar{\mathbb Q},PSL_2\RR\times L^+)/Cont(\hat{G},PSL_2\RR ),$$
where $\alpha :\hat{G}\to PSL_2\RR$ acts on $\tilde\rho$ by
$$(\alpha \tilde\rho)(\gamma,t,q)=\bigl(\alpha(t\gamma^{-1})
\circ\rho(\gamma ,t)\circ\alpha^{-1}(t)\bigr)\times \bigl(\alpha
(t)h(t,q)\bigr).$$ The topology on $\tilde T(\S_G )$ is the
quotient of the topology on $Hom(G\times\hat{G}\times\bar{\mathbb
Q},PSL_2\RR\times L^+)$. The following is immediate:

\vskip .2in

\noindent {\bf Proposition 5.2}~~\it Forgetting decoration induces
a continuous surjection $\tilde {T}(\S )\to T(\S )$.\rm

\vskip .2in

\noindent

Given a geodesic in ${\mathbb D}$ together with a horocycle
centered at each end of the geodesic, the intersection points of
the horocycles determine a finite segment of the geodesic. To this
geometric data, we assign the square root of the double of the
exponential of the signed length of the segment (positive if the
segment lies outside both horoballs, otherwise negative) and call
it {\it lambda length}. (For more details see the appendices.)  Given
$\tilde{\rho}$, to any $e\times t$, for $e\in\tau_{*}$ and
$t\in\hat{G}$, we can assign the lambda length of the image geodesic
$\phi_t(e)$ and corresponding horocycles $h(t,q_1)$ and
$h(t,q_2)$, where $q_1,q_2$ are images of the endpoints of $e$
under $\phi_t$. Thus, we obtain a natural mapping
$\lambda:\tilde T(\S )\to({{\mathbb R}}_{>0}^{\tau _*})^{\hat G}$
which assigns to a function
$\tilde\rho:G\times\hat{G}\times\bar{\mathbb Q}\to PSL_2\RR\times
L^+$ the lambda lengths corresponding to the edges $(e,t)$, for
$e\in\tau_{*}$ and $t\in\hat{G}$.

We consider ${{\mathbb R}}_{>0}^{\tau _*}$ with the strong
topology induced by taking the supremum over the coordinates. Let
$Cont(\hat{G},{{\mathbb R}}_{>0}^{\tau _*})$ be the space of all
continuous mappings with the compact-open topology. Denote by
$Cont^{G}(\hat{G},{{\mathbb R}}_{>0}^{\tau _*})$ the subset of
$Cont(\hat{G},{{\mathbb R}}_{>0}^{\tau _*})$ which consists of all
elements that are invariant with respect to the action of $G$. In
other words, $f\in Cont^G(\hat{G},{{\mathbb R}}_{>0}^{\tau _*})$ if
$f\in Cont(\hat{G},{{\mathbb R}}_{>0}^{\tau _*})$ and
$$
f(t\gamma^{-1},\gamma (e))=f(t,e)
$$
for each $\gamma\in G$ and $e\in\tau^{*}$.

\vskip .2in

\noindent{\bf Theorem 5.3}~~\it The assignment of lambda lengths
is a homeomorphism onto
$$\lambda :\tilde{T}(\S_G )\to Cont^G(\hat{G}, {{\mathbb R}}_{>0}^{\tau
_*});$$ that is, we obtain a parametrization of $\tilde{T}(\S_G)$.\rm

\vskip .2in

\noindent{\bf Proof}~~We first show that $\lambda$ is surjective.
If $f\in Cont^G(\hat{G},{{\mathbb R}}_{>0}^{\tau _*})$, then we
must produce $\tilde{\rho}=\rho\times h$ such that $\lambda
(\tilde{\rho} )=f$. Let $\gamma_1$ and $\gamma_2$ be fixed
generating hyperbolic elements in the once punctured torus group $G$
that map $e_{1/2}$ and $e_{2/1}$ onto $e_{-2/1}$ and $e_{-1/2}$,
respectively.

By Lemma A.1(i), given $f(t,e_{1/2})$, $f(t,e_{2/1})$ and
$f(t,e_0)$ there exists a unique choice of horocycles $h(t,-1)$,
$h(t,-i)$ and $h(t,1)$ based at $(-1,t)$, $(-i,t)$ and $(1,t)$ in
${\mathbb D}\times\hat{G}$ such that induced lambda lengths on
$e_{1/2}\times t$, $e_{2/1}\times t$ and $e_0\times t$ are equal
to the above values of $f$. Further by Lemma A.1(i), given $f(e_{-1/2},t)$ and
$f(e_{-2/1},t)$ there exists a unique horocycle $h(t,i)$ in
${\mathbb D}\times\hat{G}$ which induces lambda lengths
$f(t,e_{-1/2})$ and $f(t,e_{-2/1})$ on $e_{-1/2}\times t$ and
$e_{-2/1}\times t$, respectively. We continue this process
inductively and obtain $h(t,e)$ for all $e\in\tau_{*}$ and
$t\in\hat{G}$ such that whenever $p,q\in\bar{\mathbb{Q}}$ are
endpoints of an edge $e\in\tau_{*}$ then the corresponding lambda
length for $h(t,p)$ and $h(t,q)$ equals $f(t,e)$.

We define $\phi_t|_{\bar{\mathbb{Q}}}$ by mapping $(p,t)$ onto the
center of $h(t,p)$, for each $p\in\bar{\mathbb{Q}}$. We show that
$\phi_t$ extends to a quasiconformal map $\phi_t
:{\mathbb D}\times\hat{G}\to{\mathbb D}\times\hat{G}$ which is continuous
in
$t\in\hat{G}$. By the invariance under $G$, $f$ descends to a
continuous function on the compact set $(\tau_{*}\times\hat{G})/G$.
Therefore, the lambda lengths on each ${\mathbb D}\times\{ t\}$ defined by
$f$ are pinched (cf. Lemma 6.1 and Appendix A). A theorem of Penner
and Sullivan \cite{P1} (cf. Theorem~A.2) guarantees that $\phi_t$
extends to a quasiconformal map. The method of \cite{P1} shows
that $\phi_t$ varies continuously in $t$ upon replacing $\phi_t$
with the barycentric extension of its boundary values.

We define $\rho (\gamma ,t )$ to be the unique element of
$PSL_2\RR$ which maps $(-1,t)$, $(-i,t)$ and $(1,t)$ onto
$\phi_{t\gamma^{-1}}(\gamma (-1),t\gamma^{-1})$,
$\phi_{t\gamma^{-1}}(\gamma (-i),t\gamma^{-1})$ and
$\phi_{t\gamma^{-1}}(\gamma (1),t\gamma^{-1})$, where we identify
${\mathbb D}\times\{ t\}\equiv {\mathbb D}\times\{
t\gamma^{-1}\}\equiv{\mathbb D}$ for
$\gamma\in\hat{G}$. We need to show that $\rho (\gamma ,t
)=\phi_{t\gamma^{-1}}\circ\gamma\circ\phi_t^{-1}$ for all
$t\in\hat{G}$, and it is enough to show that the equality holds on the
edges of $\phi_t(\tau_{*})$. Any $e\in\phi_t(\tau_{*})$ is
obtained by forming a unique chain of edges from the base triangle with
vertices $1$, $-1$ and $-i$ according to the values of $f$ on the
edges of $\tau_{*}$. In the same way, we obtain the edge
$\phi_{t\gamma^{-1}}\circ\gamma\circ\phi_t^{-1} (e)$ on ${\mathbb D}\times
t\gamma^{-1}$. Since the function $f$ is invariant under the
action of $G$, it follows that the values of $f$ on the chain of
edges from the base triangle to $e$ on ${\mathbb D}\times t$ are equal to
the values of $f$ on the corresponding edges of the chain from the
image of the base triangle under
$\phi_{t\gamma^{-1}}\circ\gamma\circ\phi_t^{-1}$ to the edge
$\phi_{t\gamma^{-1}}\circ\gamma\circ\phi_t^{-1}(e)$ on ${\mathbb D}\times
t\gamma^{-1}$. On the other hand, $\rho (\gamma ,t )\in
PSL_2\mathbb{R}$ preserves the geometric construction using
the horocycles and it agrees with
$\phi_{t\gamma^{-1}}\circ\gamma\circ\phi_t^{-1}$ on the base
triangle. Thus, $\rho (\gamma ,t)$ and
$\phi_{t\gamma^{-1}}\circ\gamma\circ\phi_t^{-1}$ agree on $e$,
hence $\rho$ is a representation which automatically satisfies
Property 2. Property 3 is also established. Finally, the
representation $\rho$ is continuous because the function $f$ is
continuous, and Property 1 follows.

Since $f$ is pinched, it follows from \cite{P1} (cf. Theorem~A.2) that
$h(t,\bar{\mathbb Q})$ is discrete and radially dense. The continuity of
$h(\cdot , q):\hat{G}\to L^{+}$ follows by the continuity of $f$ and by
continuity of the formulas in Lemma~A.1(i). Property 6 follows from
invariance of lambda lengths and invariance of $f$.

It follows that  $\lambda (\tilde{\rho})=f$, and hence $\lambda$ is
surjective. If
$f$ and $f_1$ are close, then the above construction yields
$\tilde{\rho}$ and $\tilde{\rho}_1$ which are close in
$\tilde{T}(\S_G)$. This holds by the continuity of the
construction.

The map $\lambda$ is continuous because lambda lengths are
invariant under $G$ and depend continuously on the representation
$\tilde{\rho}=\rho\times h$.

It remains to show that $\lambda$ is injective. A representation
$\tilde{\rho}$ up to conjugation by $\alpha\in
Cont(\hat{G},PSL_2\RR )$ is determined uniquely by the decorations
$h(t,-1)$, $h(t,-i)$, $h(t,1)$ and $h(t,i)$, for $t\in\hat{G}$. On
the other hand, the above decorations and invariance with respect
to $G$ uniquely determine the function $f=\lambda (\tilde{\rho})\in
Cont^G(\hat{G},\RR_{>0}^{\tau^{*}})$, whence $\lambda$ is indeed
one-to-one. $\Box$

\vskip .2in

\noindent{\bf Remark 5.4}~One may parametrize other transverse
structures on the solenoid in analogy, where conditions other than
continuity are imposed on the  ``lambda functions'' $\hat G\to
{{\mathbb R}}_{>0}^{\tau _*}$.

\vskip .2in

The above parametrization of $\tilde{T}(\S_G )$ immediately
implies density of TLC decorated structures on  $\S$ since
$\hat{G}$ is a Cantor set in which $G$ is dense. This is in
parallel to Theorem 3.8:

\vskip .3 cm

\paragraph{\bf Corollary 5.5} {\it The union of the natural lifts of
the decorated Teichm\"uller spaces of all finite punctured
surfaces is dense in the decorated Teichm\"uller space
$\tilde{T}(\S_G )$.}

\section{Convex hull construction for the solenoid}

\vskip .2in

In fact, the results of Appendix~A automatically apply
to any continuous lambda length function because of the following:

\vskip .2in

\paragraph{\bf Lemma 6.1}{\it Continuity of a $G$-invariant $\lambda:\hat
G\to{{\mathbb R}}_{>0}^{\tau _*}$ implies that $\lambda _t:\tau
_*\to{{\mathbb R}}_{>0}$ is pinched for each $t\in \hat G$.}

\vskip .2in

\paragraph{\bf Proof}~Continuity of $\lambda$ means that $\forall
K~\exists N~\forall e\in\tau _*~\forall\gamma\in G_N$, we have
$$1+K^{-1}\leq {{\lambda _t(e)}\over{\lambda _t(\gamma e)}}\leq 1+K.$$
Take $K=1/2$ and its corresponding $N$.  A fundamental domain for
$G_N$ has only a finite collection of lambda lengths, and any
other is at most three halves and at least one half times a lambda
length in this finite set.$\Box$

\vskip .2in

Thus, if $\lambda:\hat G\to{{\mathbb R}}_{>0}^{\tau _*}$ is
continuous, then each $\lambda _t:\tau _*\to{\mathbb R}_{>0}$ is pinched, so in the notation
of Appendix~A, the corresponding $h_t:\tau _*\to L^+$ has drd
image ${\mathcal B}_t$ and closed convex hull $C_t$, which
projects to the cell decomposition $\tau _t$ of ${{\mathbb D}}$.
Presumably, cells in the decomposition could be
infinite sided. Furthermore, the characteristic map interpolates a
quasi-symmetric mapping $\phi _t:S^1\to S^1$.

Given a tesselation $\tau$ of ${{\mathbb D}}$ with $\tau ^\infty={\mathbb
Q}\cup\{\infty\}=\bar{\mathbb Q}=\tau _*^\infty$, define
\begin{equation*}
\begin{array}l
\stackrel{\circ}{\mathcal C}(\tau )=\{\lambda _t\in \tilde
T({\mathcal H}): \tau _t=\phi _t(\tau)\} \\ \ \ \cap\\{\mathcal
C}(\tau )=\{\lambda _t\in \tilde T({\mathcal H}): \tau _t\subseteq\phi
_t(\tau)\}\\ \ \ \cap\\ \tilde T({\mathcal H})\approx Cont ^G(\hat
G,{{\mathbb R}}_{>0}^{\tau _*}).
\end{array}
\end{equation*}

Furthermore, define the {\it classical locus} ${\mathcal
L}\subseteq \tilde T({\mathcal H})$ as the subspace consisting of
all TLC structures on ${\mathcal H}$, or equivalently, as the
union over all finite-index subgroups $K < G$ of the image of
$\tilde T({{\mathbb D}}/K)$ in $\tilde T({\mathcal H})$.  The
corresponding subspace ${\mathcal L}\subseteq Cont ^G(\hat
G,{{\mathbb R}}_{>0}^{\tau _*})$ is described by the collection of
all TLC lambda length functions $\lambda _t:\tau _*\to {{\mathbb
R}}_{>0}$ where there is some finite-index $K <G$ so that
$\lambda _t$ is $K$-invariant, for all $t\in\hat{G}$.  Define
a {\it TLC tesselation} to be a tesselation $\tau$ of ${{\mathbb D}}$ that
is invariant under some finite-index subgroup $K<G$ with ideal points  $\tau ^\infty
=\bar{\mathbb Q}$.

For any fixed TLC tesselation $\tau$, a point of ${\mathcal
C}(\tau )-\stackrel{\circ}{\mathcal C}(\tau )$ corresponds to
$\lambda _t:\tau _*\to{{\mathbb R}}_{>0}$ whose convex hull $C_t$
has at least one face which is not triangular, i.e., at least four
points of ${\mathcal B}_t$ are coplanar; thus, $\tau _t$ is not
a tesselation, but rather a paving of ${{\mathbb D}}$. Given an arbitrary
point in $\tilde{T}(\S_G )$, it is presumably possible that it does not
belong to ${\mathcal C}(\tau )$ for any TLC tesselation $\tau$.  Notice,
however, that for any point of
$\tilde{T}(\S_G )$, we have
$\tau _t^\infty\subseteq\bar{\mathbb Q}$. On the other hand, we
obtain a generically simple picture:

\vskip .2in

\paragraph{\bf Theorem 6.2}{\it The subspace $\stackrel{\circ}
{\mathcal C}(\tau )$ of $\tilde T({\mathcal H})$ is open for each
TLC tesselation $\tau$, and
 $\cup _\tau~\stackrel{\circ} {\mathcal
C}(\tau )$ is a dense open subset of $\tilde T({\mathcal H})$,
where the union is over all TLC tesselations $\tau$ of ${\rm D}$.}

\vskip .2in

\paragraph{\bf Proof}~Fix $\lambda \in Cont ^G(\hat G,{{\mathbb
R}}_{>0}^{\tau _*})$ and $\epsilon >0$. There exists a TLC
$\lambda^{'}\in Cont ^G(\hat G,{{\mathbb R}}_{>0}^{\tau _*})$, say
with corresponding finite-index subgroup $K <G$, so that
$\|\lambda -\lambda^{'}\|_{\infty} <\epsilon /2$ by Corollary 5.5.
By the classical theory \cite{P2} working on the surface
${{\mathbb D}}/K$, there exists $\lambda''$ which belongs to
$\stackrel{\circ}{\mathcal C}(\tau )$ for some TLC tesselation
$\tau$ and which satisfies $\|\lambda'' -\lambda^{'}\|_{\infty}
<\epsilon /2$. Therefore, the density follows.

Turning to the proof that each $\stackrel{\circ} {\mathcal C}(\tau
)$ is open and in order to distinguish coordinates, let us now
denote
$$Cont_\lambda (\hat G,{{\mathbb R}}_{>0}^{\tau _*})
=Cont ^G (\hat G,{{\mathbb R}}_{>0}^{\tau _*})\approx \tilde T({\mathcal
H}).$$ Let $Cont_\sigma (\hat G,{{\mathbb R}}^{\tau _*})$ denote the
abstract space of continuous $G$-equivariant ${{\mathbb R}}$-valued
function, and define
\begin{equation*}
\begin{array}l
\Phi:Cont_\lambda (\hat{G},{{\mathbb R}}_{>0}^{\tau _*}) \to Cont_\sigma
(\hat{G},{{\mathbb R}}^{\tau })\\
\hskip 2.6 cm \lambda _t \mapsto \sigma _t,
\end{array}
\end{equation*}
where the ``simplicial coordinate function'' $\sigma _t$ is
defined for each $t\in\hat G$ in accordance with the formula in
Lemma A.1(v).  By definition of simplicial coordinates and the convex
hull construction, in fact we have

\begin{equation*}
\begin{array}l
\stackrel{\circ}{\mathcal C}(\tau )= \Phi^{-1} ( Cont_\sigma
(\hat G,{{\mathbb R}}_{>0}^{\tau })),\\
{\mathcal C}(\tau )= \Phi ^{-1} (Cont_\sigma (\hat G,{{\mathbb R}}_{\geq
0}^{\tau })),
\end{array}
\end{equation*}
so the openness assertion of the proposition thus follows from
continuity of $\Phi$. $\Box$

\vskip .2in

\noindent {\bf Remark 6.3} In the classical case \cite{P2} (cf.
Appendix~B), $\tau$ is a paving of ${{\mathbb D}}$ arising as the lift of
a paving of some punctured surface ${{\mathbb D}}/\Gamma$.  The
intersection ${\mathcal L}
\cap {\mathcal C}(\tau )$ is mapped homeomorphically by $\Phi$
onto the set of TLC elements of $Cont_\sigma (\hat G,{{\mathbb
R}}_{\geq 0}^{\tau })$ so that there are no (finite) cycles of
triangles $t_i,\ldots t_{i+k}$ with $t_i$ and $t_{i+k}$ in the
same $\Gamma$-orbit and with vanishing simplicial coordinates on
each $t_i,\ldots ,t_{i+k-1}$.  Each component of ${\mathcal
L}\cap {\mathcal C}(\tau )$ is homeomorphic to an open simplex
together with certain of its faces.

\vskip .2in

\section{The Modular Group}

We define the modular group of the punctured solenoid in analogy to
the modular group of the compact universal hyperbolic solenoid
\cite{Odd},\cite{Od}, \cite{MS}.

\vskip .3 cm

\paragraph{\bf Definition 7.1}
The {\it modular group} $Mod(\S_G )$ of the punctured solenoid
$\S_G$ consists of quasiconformal self-maps of $\S_G$ modulo
isotopies which are bounded in the hyperbolic metric on leaves.
The {\it baseleaf preserving modular group} $Mod_{BLP}(\S_G )$ is
a subgroup of $Mod(\S_G )$ which consists of all isotopy classes
of quasiconformal self-maps of $\S_G$ which fix the baseleaf.

\vskip .3 cm

Note that a quasiconformal self-map $f$ of $\S_G$ necessarily sends
punctures onto punctures on $\S_G$. Indeed, since an arbitrary
quasiconformal map of the hyperbolic plane onto the hyperbolic
plane is a quasi-isometry \cite{EMM}, a geodesic on a leaf which
ends in the puncture is mapped by $f$ onto a quasigeodesic on the
image leaf. If the endpoint of the quasigeodesic on the image leaf
is not a puncture then it returns to a compact set of $\S_G$
infinitely often, and this contradicts that $f$ is proper.

C. Odden \cite{Odd}, \cite{Od} showed that the baseleaf preserving modular
group of the compact universal hyperbolic solenoid is isomorphic
to the group of virtual automorphisms of the fundamental group of
a compact surface with genus greater than one. We use his method
to prove a similar statement for $\S_G$.

\vskip .3 cm

\paragraph{\bf Theorem 7.2}{\it The restriction of $Mod_{BLP}(\S_G )$ to
the baseleaf is isomorphic to the virtual automorphism group of
$G$.}

\vskip .1 cm

\paragraph{\bf Proof} We briefly discuss the extension of Odden's argument
(see \cite[Theorem 4.6]{Od}) to the case of the punctured
solenoid. A quasiconformal self-map $f$ of the compact universal
solenoid is uniformly continuous, and the compact solenoid admits
a finite covering by $\epsilon$ balls in the topology coming from
the representation ${{\mathbb D}}\times\hat{G}/G$, for a cocompact
Fuchsian
$G$.  These are two crucial ingredients in Odden's proof. We
replace compactness of the solenoid by considering a compact
subset $\tilde{X}$ of the punctured solenoid obtained by lifting
of the complement $X\subset{\mathbb D} /G$ of a horoball neighborhood of
the cusp on ${{\mathbb D}} /G$, where $G$ is the punctured torus group.
Further, $\tilde{X}_1=f(\tilde{X})$ is a compact subset of $\S_G$
which is contained in another compact set $\tilde{X}_2$ of the
same kind as $\tilde{X}$. In this case, $f:\tilde{X}\to
\tilde{X}_1$ is uniformly continuous and both
$\tilde{X},\tilde{X}_1$ can be covered by finitely many $\epsilon$
balls as required, and they furthermore admit a local product
structure similar to the compact solenoid. Moreover, each closed
curve based at a point in $ X\subset{{\mathbb D}} /G$ can be homotoped
into a closed curve in $ X$ itself. These details and observations
allow the further application of Odden's proof to the current
setting. $\Box$

\vskip .3 cm

For a compact solenoid, the isotropy group of a point in the
Teichm\"uller space in the baseleaf preserving modular group is
infinite (see \cite{Odd}, \cite{Od} for TLC
hyperbolic solenoids and see \cite{MS} for non-TLC hyperbolic
solenoids). Moreover, each isotropy group contains a subgroup
isomorphic to the surface group. There is a countable set of TLC
hyperbolic solenoids for which the isotropy group is isomorphic to
a dense subgroup of $PSL_2(\mathbb{R})$ (see \cite{Odd},\cite{Od}). These
statements hold for $T(\S_G )$ as well by obvious generalizations of the
proofs.

\vskip .3 cm

We investigate the isotropy subgroups of  $Mod_{BLP}(\S_G
)$ for a point $\tilde{\S}\in\stackrel{\circ}{\mathcal C}(\tau )\subseteq\tilde{T}(\S_G )$
for $\tau$ a TLC
tesselation. Note that the pull-back $\tau _t$ under the
marking map $\phi_t$ of the extreme edges of the boundary of the
convex hull $C_t$ in Minkowski space for a decorated hyperbolic
punctured solenoid $\tilde{\S}$ is invariant under the action of
$Mod_{BLP}(\S_G )$ by definition. That is, an element $h\in Mod_{BLP}(\S_G )$
which fixes $\tilde{\S}\in\stackrel{\circ}{\mathcal C}(\tau )$ must map $\tau _t$
onto itself. Since the projection $\pi :\tilde{T}(\S_G )\to T(\S_G
)$ commutes with the action of $h$, it follows that $h$ fixes $\S
=\pi (\tilde{\S})\in T(\S_G )$, and so $h$ is an isometry of $\S$.

>From now on, we restrict the action to the baseleaf, which is
identified with the unit disk ${{\mathbb D}}$. A homeomorphism
$h:S^1\to S^1$ extends diagonally to a map of the space of
geodesics $\{ x\times y\in S^1\times S^1: x\neq y\}$. To study the
isotropy groups, we need the following simple but important
rigidity property of Farey tesselation:

\vskip .4 cm

\paragraph{\bf Lemma 7.3}\it
Assume that $h:S^1\to S^1$ is a homeomorphism which maps the Farey
tesselation $\tau_{*}$ onto itself. Then $h$ is an element of
$PSL_2(\mathbb{Z})$.\rm

\vskip .4 cm

\paragraph{\bf Proof} The
homeomorphism $h:S^1\to S^1$ maps $\bar{\mathbb{Q}}$ onto itself
because the endpoints of any geodesic in $\tau_{*}$ are mapped
onto the endpoints of another geodesic in $\tau_{*}$. We show that
$h$ must belong to $PSL_2\mathbb{Z}$. Fix $a\in\tau_{*}$ and let
$T\subset {{\mathbb D}}-\tau_{*}$ be an ideal triangle with $a$ in its
frontier. Let $\gamma$ be the unique element of $PSL_2\mathbb{Z}$ such
that $\gamma (a)=h(a)$ and $\gamma (T)\cap h(T)\neq\emptyset$. If
$\gamma\neq h$ then there exists $a\in\tau_{*}$ and a triangle $T$ complementary to $\tau _*$
whose
frontier contains $a$ such that $\gamma^{-1}\circ h(T)\neq T$. In
this case, $\gamma^{-1}\circ h$ maps an edge of $T$ which is in
$\tau_{*}$ onto a geodesic not in $\tau_{*}$. However, both
$\gamma$ and $h$ preserve $\tau_{*}$ which gives a contradiction.
$\Box$

\vskip .2 cm

Fix a{\it TLC tesselation} $\tau$ of ${{\mathbb D}}$, i.e. $\tau$
is an ideal triangulation of ${\mathbb D}$ which is invariant
under a finite index subgroup of $PSL_2\mathbb{Z}$ and has ideal
vertices $\tau^{\infty}=\bar{\mathbb{Q}}$. Denote by $Aut(\tau )$
the subgroup of the baseleaf preserving modular group
$Mod_{BLP}(\S_G )$ which fixes $\tau$ setwise. In particular,
$Aut(\tau_{*})=PSL_2\mathbb{Z}$ (by above lemma) acts simply
transitively on the oriented edges of $\tau _*$ and fixes the
basepoint $\S_G\in T(\S_G )$ and the basepoint
$\tilde{\S}_G\in\tilde{T}(\S_G )$, where the decoration on
$\tilde{\S}_G$ is given by assigning lambda length unity to each
edge of $\tau_{*}$.

Note that if $\tau$ is a TLC tesselation invariant under some group $K$, say without elliptics, then
$\tau/K$ is a well-defined tesselation of the surface $S={\mathbb D}/K$.  Furthermore, if
$f$ is an isotopy class of  homeomorphisms of $S$ fixing $\tau/K$ setwise, then a lift of $f$ to the universal cover $\mathbb D$ fixes the pre-image of $\tau /K$ in ${\mathbb D}$.  It follows that every automorphism group of every TLC tesselation invariant under a finite index subgroup $K<PSL_2({\mathbb Z})$ lies as a subgroup of $PSL_2({\mathbb Z})$.  (The same is not true for more general pavings of $S$ as one can easily check directly using simplicial coordinates, cf. Lemma~A.1(v).)

We define the {\it characteristic map} $h=h(\tau ,e)$ for a
tesselation $\tau$ with a distinguished oriented edge and recall
that $\tau ^\infty\subset S^1$ denotes the set of ideal points of
the tesselation $\tau$. Define $h$ to map the initial and final
points of the {\it DOE} of $\tau_{*}$ onto the initial and final
points of the {\it DOE} of $\tau$. Each {\it DOE} is the common
boundary of two complementary triangles, one to the left and one
to the right with respect to the orientations on the {\it DOES}.
Map the third vertices of the respective triangles in ${\mathbb D}
-\tau_{*}^{\infty}$ to the third endpoints of corresponding
triangles in ${\mathbb D} -\tau$.   Continue in this manner
mapping third points of triangles to recursively define $h:\tau
_*^\infty \to \tau^\infty$. By construction, $h$ is monotone and
hence interpolates a homeomorphism $h$ of $S^1$. (See Appendix A
for more details.)

\vskip .1 cm

There is the following immediate corollary of the previous result:

\vskip .3cm

\noindent{\bf Corollary 7.4}~\it Suppose that $\tau$ is a TLC
tesselation of ${\mathbb D}$ with distinguished oriented edge $e$ and
corresponding characteristic map $h=h{(\tau
,e)}:S^1\to S^1$.  Then $Aut(\tau )=h\circ PSL_2{\mathbb Z}\circ
h^{-1}$. Furthermore, if $\phi=h\circ\gamma\circ h^{-1}$, for
$\gamma\in PSL_2{\mathbb Z}$, then $h{(\tau ,\phi (e))}=h(\tau
,e)\circ \gamma$.\rm

\vskip .3 cm

\paragraph{\bf Lemma 7.5} {\it Let $\tau$ be a
TLC tesselation of ${{\mathbb D}}$, i.e., $\tau$ is invariant under a
finite-index subgroup $K$ of $PSL_2\mathbb{Z}$ with ideal points
$\tau^{\infty}=\bar{\mathbb{Q}}$. Then the characteristic map of
$\tau$ conjugates a finite-index subgroup $H$ of $PSL_2\mathbb{Z}$
onto $K$.}

\vskip .1 cm

\paragraph{\bf Proof} Let $K^{'}$ be the finite-index subgroup
of $K$ which contains no elliptic elements. Choose a fundamental
polygon $P'$ for $K^{'}$ whose boundary consists of edges in
$\tau$. Thus,  the interior of $P'$ is divided into ideal
triangles by the edges of $\tau$. To describe subgroup $H$,
consider a combinatorially equivalent polygon $P$ which is
comprised of ideal triangles among ${{\mathbb D}}-\tau_{*}$. The
combinatorial correspondence between $P'$ and $P$ gives the
pairing of the boundary of $P$ using the pairing of $P'$ under
$K^{'}$. For each pair $a\times b\in \tau _*\times \tau _*$ on the
boundary of $P$, there exists a unique $\gamma_{a,b}\in
PSL_2\mathbb{Z}$ such that $\gamma_{a,b}(a)=b$ and
$\gamma_{a,b}(P')\cap P'=\{ b\}$. By Poincar\'e's fundamental
polygon theorem, the group $H$ generated by all side pairings
$\gamma_{a,b}$ of the corresponding boundary edges of $P$ is
Fuchsian with a fundamental polygon equal to $P$. Thus, the
Riemann surfaces ${{\mathbb D}}/K^{'}$ and ${{\mathbb D}}/H$ are
homeomorphic, and we may choose a quasiconformal map between them which
maps edges in $P'$ onto the edges of $P$. This map lifts to ${{\mathbb
D}}$, and its restriction to $S^1$ is the desired quasisymmetric
conjugation $h$ between $K^{'}$ and $H$. The map $h$ differs from the
characteristic map by a pre-composition with an element of
$PSL_2\mathbb{Z}$.

Now, if $K'\neq K$, then we can choose $P'$ such that the orbit of
$P'$ under $K^{'}$ is invariant under the full group $K$. This is
equivalent to requiring that a fundamental set for the action of
$K$ is contained in $P'$ and finitely many of its translates under
some elliptic elements of $K$ cover $P'$. Any elliptic element
$\delta\in K-K^{'}$ is conjugated by $h$ to a finite order
quasisymmetric element $h\circ\delta\circ h^{-1}$ of ${{\mathbb D}}$ which
preserves $\tau_{*}$, whence it must be an element of
$PSL_2\mathbb{Z}$ by Lemma 7.3.
 $\Box$

\vskip .3 cm

>From the above we obtain

\vskip .1 cm

\paragraph{\bf Theorem 7.6} {\it $Mod_{BLP}(\S_G )$ acts
transitively on $\{ \stackrel{\circ}{\mathcal C}(\tau):\tau ~{\rm
is}~TLC\}$.  Furthermore, the isotropy subgroup in $Mod_{BLP}(\S_G
)$ of $\stackrel{\circ}{\mathcal C}(\tau )$ is isomorphic to
$Aut(\tau )$ and is moreover quasi-conformally conjugate to
$PSL_2({\mathbb Z})$.}

\vskip .1 cm

\paragraph{\bf Proof} In Lemma 7.5, we showed that a
TLC tesselation $\tau$ of ${{\mathbb D}}$ which is invariant under a
finite subgroup $K$ of $PSL_2\mathbb{Z}$ can be mapped to $\tau_{*}$ by a
quasisymmetric map $h$ which conjugates $K$ onto a finite-index
subgroup of $PSL_2\mathbb{Z}$. Thus, $h$ is a virtual
automorphism, and by Theorem 7.2, it defines an element of
$Mod_{BLP}(\S_G )$ which sends the cell corresponding to
$\tau_{*}$ onto the cell corresponding to $\tau$. Transitivity
follows.

The identification of the isotropy group of
$\stackrel{\circ}{\mathcal C}(\tau )$ with $Aut(\tau )$ is induced
by identifying both these groups with the isometry group of the
point of $\tilde T(\S _G)$ described via Theorem 5.3 with all
lambda lengths constant equal to unity.  This point lies in
$\stackrel{\circ}{\mathcal C}(\tau )$ (as one checks with
simplicial coordinates), and has combinatorial symmetry group
given by $Aut(\tau )$ and decorated hyperbolic (or conformal) symmetry group given
by the isotropy subgroup.  The last part then follows from
Corollary 7.4. $\Box$

\vskip .3cm

We show that the space $\tilde{T}(\S_G ) /Mod_{BLP}(\S_G )$ is
``essentially'' Hausdorff, more precisely, at least an open dense
subset is Hausdorff. This implies that no orbit under
$Mod_{BLP}(\S_G )$ is dense. The analogue of the Ehrenpreis
conjecture is thus very false for $\tilde{T}(\S_G )$.

\vskip .3 cm

\paragraph{\bf Theorem 7.7} {\it The quotient $\cup_{\tau}
\stackrel{\circ}{\mathcal C}(\tau ) /Mod_P(\S_G )$ is Hausdorff,
where the union is over all TLC tesselations $\tau$. Moreover, no
orbit of a point in $\tilde{T}(\S_G )$ is dense.}

\vskip .1 cm

\paragraph{\bf Proof}
Fix an identification of the baseleaf of $\S_G$ with the unit disk
${\mathbb D}$. Since the points in $\tilde{T}(\S_G )$ are defined up to
post composition by a conformal map, we may identify
the image of the baseleaf under $f:\S_G\to\tilde{\S}$ with ${\mathbb D}$
such that $f$ pointwise fixes each of $\pm 1,-i$. Thus, $\tilde{T}(\S_G )/Mod_P(\S_G
)$ is mapped to the space ${\mathcal B_{drd}}$ of all discrete, radially
dense countable subsets of the light cone
$L^{+}$ in Minkowski three space containing three
points which project to $1,-1,-i$ on $S^1$. The map
$\tilde{T}(\S_G )/Mod_P(\S_G )\to {\mathcal B}_{drd}$ is
continuous for the quotient topology on $\tilde{T}(\S_G
)/Mod_P(\S_G )$ and the Hausdorff topology on ${\mathcal
B}_{drd}$.

We show that the restriction of the map to $\cup_{\tau}
\stackrel{\circ}{\mathcal C}(\tau ) /Mod_{BLP}(\S_G )$ is
injective which immediately implies that the space is Hausdorff.
Let $f\in ~\stackrel{\circ}{\mathcal C}(\tau )$ and $f_1\in
~\stackrel{\circ}{\mathcal C}(\tau_1 )$ have the
same image $\beta$ in ${\mathcal B}_{drd}$. The convex hull
construction for $\beta$ yields an ideal triangulation of
${\mathbb D}$ which pulls back to $\tau$ and $\tau_1$ by the maps $f$ and
$f_1$, respectively. Since $\tau$ and $\tau_1$ are TLC
tesselation, there exist $h,h_1\in Mod_{BLP}(\S_G )$ such that
$h(\tau_{*} )=\tau$ and $h_1(\tau_{*} )=\tau_1$ by Theorem 7.6.
Thus, there exists $\gamma\in PSL_2\mathbb{Z}$ such that $f\circ
h\circ\gamma=f_1\circ h_1$, by Lemma 7.3.  Finally, $f$ and $f_1$
are therefore in the same orbit of $Mod_{BLP}(\S_G )$, and the injectivity
follows.

It remains to show that no orbit is dense. An orbit of a point in
$\stackrel{\circ}{\mathcal C}(\tau ) $ is not dense in
$\tilde{T}(\S_G )$ because it is not dense in
$\stackrel{\circ}{\mathcal C}(\tau ) $ by the above. The orbit of a
point outside $\cup_{\tau} \stackrel{\circ}{\mathcal C}(\tau )$ is not dense
because it does not meet this open dense set by modular invariance of the convex
hull construction.  $\Box$

\vskip .2cm

\section{Generators of the Modular Group}

To describe generators of $Mod_{BLP}({\mathcal H}_G)$, we shall require
certain elementary moves on TLC tesselations as follows.

\vskip .3 cm

\paragraph{\bf Definition 8.1} Let $K$ be a
finite-index subgroup of $G$, let $\tau$ be a tesselation of
${{\mathbb D}}$ which is invariant under $K$ with
ideal points $\tau^{\infty}=\bar{\mathbb{Q}}$, and suppose that $e$ is a fixed
unoriented edge of $\tau $ with distinct endpoints. Define a new tesselation $\tau '$ as
follows: for each $f\in\tau -Ke$, there is an identical edge
$f\in\tau '$; for each $f\in Ke$, consider the quadrilateral $P$
with diagonal $f$ comprised of the triangles on either side of $f$
complementary to $\tau $, and let $f'$ denote the other diagonal
of $P$; for each edge $f\in Ke$, there is a corresponding dual
edge $f'\in\tau '$. The resulting tesselation $\tau '$ is clearly
also invariant under $K$.  We say that $\tau '$ arises from
$\tau $ by the {\it ($K$-)equivariant Whitehead move along
$e\in\tau $}.

Furthermore, if $\tau $ is a tesselation with {\it DOE} $d\in\tau
$, then we may induce a {\it DOE} $d'$ on $\tau '$ as follows: if
$d\notin Ke$, then $d'=d$ as oriented edges, while if $d\in Ke$,
then there is the unique orientation on $d'$ so that correctly
oriented tangent vectors to $d,d'$ give a positive basis for the
tangent space at $d\cap d'$.

Thus, for any TLC tesselation $\tau $ invariant by a group $K$,
there is a corresponding equivariant Whitehead move for each edge
$e$ of $\tau$, and there is exactly one distinct equivariant move
for each $K$-orbit of edges of $\tau$.  An equivariant Whitehead
move acts not only on invariant tesselations but also on invariant
tesselations with {\it DOE}.

\vskip .3cm

\noindent {\bf Lemma 8.2}~\it Suppose that $K$ is a finite-index
subgroup of $G$ and $\tau$ is an invariant tesselation with {\it
DOE} $d$.  Perform an equivariant Whitehead move along the edge
$e$ to produce the invariant tesselation $\tau '$ with {\it DOE}
$d'$. Let $h=h(\tau ,d), h'=h(\tau ', d')$ denote the
characteristic maps and define $k=h'\circ h^{-1}$.  Then $k$ is a
quasi-symmetric map, a virtual automorphism of $PSL_2({\mathbb
Z})$, is independent of the choice of {\it DOE} $d$ on $\tau$, and
$k(\tau)=\tau '$ and $k(d)=d'$.\rm

\vskip .3 cm

\noindent{\bf Proof}~The characteristic maps $h,h'$ are
quasisymmetric by Theorem~A.1 since lambda lengths for a TLC
tesselation are pinched by Lemma 6.1, hence the composition $k$ is
quasi-symmetric as well. Likewise, $h,h'$ are virtual automorphims
by Lemma 7.5 since $\tau,\tau '$ are TLC, and hence so too is the
composition $k$.

To prove that $k$ is independent of the choice of {\it DOE}, let
$d_1$ be another {\it DOE} on $\tau$.  According to Corollary 7.4,
the characteristic maps are related by $h{(\tau,d_1)}=h (\tau
,d)\circ\gamma$, for some $\gamma\in PSL_2({\mathbb Z})$, and in
fact for the same $\gamma$, we have also $h{(\tau ',d_1')}=h{(\tau
',d')}\circ\gamma$.  Thus, $k_1=h{(\tau
',d')}\circ\gamma\circ\gamma ^{-1}\circ h^{-1}(\tau ,d)=k$ is
indeed invariant.

That $k$ maps $\tau$ to $\tau '$  and $d$ to $d'$ follow from the definition of
characteristic maps, completing the proof.
$\Box$.

\vskip .3cm

\noindent{\bf Definition 8.3}~We call $k=k(\tau ,K,e)\in
Mod_{BLP}(\S _G)$ the {\it Whitehead homeomorphism} associated
with the $K$-equivariant Whitehead move along $e\in\tau$ for any
TLC $K$-invariant tesselation $\tau$. Notice that by definition if $f\in Ke$, then $k(\tau ,K,e)=k(\tau ,K,f)$, so we regard Whitehead homeomorphisms as indexed by $K$-orbits of edges $Ke$ rather than by edges $e$.

\vskip .3cm

In contrast to the case of punctured surfaces where only certain
sequences of Whitehead moves give rise to mapping classes (namely,
the sequence must begin and end with combinatorially identical
ideal triangulations), for the punctured solenoid, each
equivariant Whitehead move {\sl does} give rise to a mapping class
$k\in Mod_{BLP}$.

As elements in the group $Mod_{BLP}(\S
_G)$, any two Whitehead homeomorphisms can be composed, but there
is a special case of geometrical significance as follows.

\vskip .3in

\noindent {\bf Definition 8.4}~ Suppose that $G>K_1>K_2$ are
nested subgroups with each finite-index in the next. Perform a
$K_1$-equivariant Whitehead move along some edge $e_1$ of the
$K_1$-invariant tesselation $\tau _1$ to get $\tau _2$.  As was
observed before, both $\tau _1$ and $\tau _2$ are $K_1$-invariant,
so in particular, $\tau _2$ is also $K_2$-invariant. Next perform
a $K_2$-equivariant Whitehead move along some edge $e_2$ of $\tau
_2$, which is a sensible geometric operation on covering spaces,
to get another TLC tesselation $\tau _3$.  In this case, we have
$k(\tau _1,K_1,e_1)\circ k(\tau _2,K_2,e_2)=h_3\circ h_1^{-1}$,
where $h_i$ is the characteristic map of $\tau _i$, for $i=1,3$,
defined with compatible {\it DOES}. We will call a composition
with this property {\it geometric}. More generally, if a finite
sequence of Whithead homeomorphisms has the same property for each
pair of consecutive terms, then it is called {\it geometric}.

\vskip .3 cm

\noindent {\bf Theorem 8.5} {\it Any element of the modular group
$Mod_{BLP}(\S_G )$ can be written as a composition $w\circ
\gamma$, where $\gamma\in PSL_2({\mathbb Z})$ and $w$ is a
geometric composition of $K$-equivariant Whitehead homeomorphisms
for some fixed $K$.   In particular, $Mod_{BLP}(\S_G)$ contains
the characteristic map of any TLC tesselation $\tau$ with {\it
DOE} and the automorphism group $Aut(\tau )$. }

\vskip .3 cm

\noindent {\bf Proof}~ We claim that the Farey tesselation
$\tau_{*}$ and an arbitrary TLC tesselation $\tau$ can be
connected by a finite geometric sequence of $K$-equivariant Whitehead moves
for some fixed $K$.  It would follow in particular that
$Mod_{BLP}(\S _G)$ and the group generated by all equivariant
Whitehead homeomorphisms have the same orbits on the set
$\{\stackrel{\circ}{\mathcal C}(\tau ):\tau~{\rm is}~TLC\}$. Thus,
for any $f\in Mod_{BLP}(\S _G)$ and any TLC tesselation $\tau$,
there is some geometric word $w$ in equivariant Whitehead homeomorphisms so
that $w^{-1}\circ f$ fixes $\stackrel{\circ}{\mathcal C}(\tau
_*)$, and hence $w^{-1}\circ f=\gamma\in PSL_2({\mathbb Z})$ by
Theorem 7.6, completing the proof of generation.

In fact, the composition of $K_1$-equivariant
and $K_2$-equivariant Whitehead moves can be equivalently
described as a finite composition of $(K_1\cap K_2)$-equivariant moves,
so by intersecting finite-index subgroups, a single group $K$ suffices as
in the statement of the theorem.

It remains to show that geometric sequences of $K$-equivariant Whitehead
moves act transitively on $K$-invariant tesselations, and this is
precisely the statement from punctured surface theory (cf. Corollary~B.6)
that sequences of Whitehead moves act transitively on the set of
all ideal triangulations of the surface ${\mathbb D}/K$.

For the proof of the last sentence, given a TLC tesselation $\tau$
with {\it DOE} and characteristic map $h$, Lemma 7.5 states that
$h$ conjugates one finite index subgroup of $PSL_2\mathbb{Z}$ onto
the other, so $h$ is indeed an element of $Mod_{BLP}(\S_G )$.
$\Box$

\vskip .3cm

\noindent{\bf Remark}~Insofar as the effect of a Whitehead move on
lambda lengths is described by a Ptolemy transformation (cf. Lemma
A.1(iii)), one derives a real-algebraic representation of
$Mod_{BLP}(\S _G)$ in analogy to \cite[Section 7]{P2}.

\vskip .3cm

We next introduce a series of relations satisfied by the
generators in Theorem 8.5, where we assume throughout that $\tau$
is a $K$-invariant tesselation of the disk for some finite-index
subgroup $K$ of $G$ (and $\tau^{\infty}=\bar{\mathbb{Q}}$):

\vskip .3cm

\noindent {\bf 1)}~[Involutivity]~~  If the
$K$-equivariant Whitehead move along
$e\in\tau$ produces $e'$ in the resulting tesselation, then $k(\tau
,K,e)\circ k(\tau,K,e')=1$.

\vskip .2cm

\noindent {\bf 2)}~[Commutativity]~~ If $e\in\tau$ and $f\in\tau$ do
not share an ideal endpoint, then $k(\tau ,K,e)\circ k(\tau ,K,f)=k(\tau
,K,f)\circ k(\tau ,K,e)$; see Figure~2.

\vskip .2cm

\noindent {\bf 3)}~[Pentagon Relation]~~If we adopt the notation for
edges in the five tesselations of the pentagon illustrated in Figure~2,
then
$$k(\tau _1,K,e_1)\circ k(\tau _2,K,f_2)\circ k(\tau _3,K,e_3)\circ
k(\tau _4,K,f_4)\circ k(\tau _5,K,e_5)=1.$$

\vskip .2cm

\noindent {\bf 4)}~[Coset Relation]~~If $H$ is a finite-index
subgroup of $K$ and $f_1,f_2,\ldots ,f_\ell$ are representatives
for the cosets of $H$ in $K$, then
$$k(\tau,K,e)=k(\tau ,H,f_1)\circ
k(\tau ,H,f_2)\circ\cdots\circ k(\tau ,H,f_\ell),$$
where the order of composition is irrelevant since the
individual Whitehead homeomorphisms commute by Relation 2).

Relation 4) holds obviously. To see that the others hold, it is
enough to show that the moves leave invariant any starting
tesselation with arbitrary {\it DOE} since Lemma 8.2 guarantees
that the choice of {\it DOE} is unimportant. Relations 2) and 3)
follow from Figure 2, say if we choose {\it DOE} to be on the
boundary, and 1) is obvious for any {\it DOE} different from the
edge $e$.

\leftskip=0ex

\vskip .3cm

\hskip .2in\includegraphics{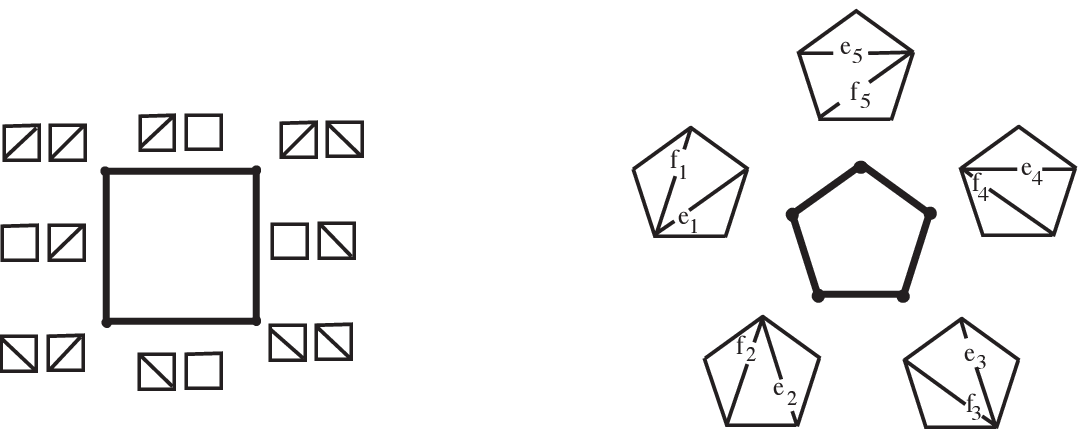}

\vskip .3cm

\hskip .5in Commutativity.\hskip 1.6in Pentagon
Relation.

\vskip .3cm

\centerline{{\bf Figure 2} Commutativity and Pentagon Relation.}

\vskip .3cm

\paragraph{\bf Definition 8.6} Let $\omega$ be a geometric composition of
equivariant Whitehead homeomorphisms which are not necessarily
equivariant with respect to the same finite index subgroup of
$PSL_2(\mathbb{Z} )$. Let $\gamma\in PSL_2(\mathbb{Z} )$. Any such
composition $\omega\circ\gamma$ is called a {\it normal form} of an element of
$Mod_{BLP}(\S_G )$.

\vskip .1 cm

We showed in Theorem 8.5 that each element of $Mod_{BLP}(\S_G )$
has a normal form, and we next address its uniqueness.

\vskip .3 cm

\noindent {\bf Theorem 8.7}\it ~~Suppose $\gamma_1 ,\gamma_2\in
PSL_2(\mathbb{Z} )$ and  $\omega_1 ,\omega_2$ are two geometric
Whitehead compositions with $\omega_1\circ\gamma_1
=\omega_2\circ\gamma_2$.  Then there is some finite index subgroup
$K$ of $PSL_2({\mathbb Z})$ with $S={\mathbb D}/K$ so that up to
Relations 1-4) for $K$, $\omega _2^{-1}\circ\omega _1=\gamma
_2^{-1}\circ \gamma _1$ is a finite composition $\phi _1\circ\phi
_2\circ\cdots\circ \phi _k$ of automorphisms $\phi _i\in Aut(\tau
_i)<PSL_2({\mathbb Z})$ of tesselations $\tau _i$ of $S$ without
{\it DOE}, for $i=1,\ldots, k$. \rm

\vskip .3cm

\noindent{\bf Proof}~Since $\omega_1\circ\gamma_1
=\omega_2\circ\gamma_2$ and $\omega_1 ,\omega_2$ are geometric, we
conclude that the image tesselations of $\tau _*$ under $\omega_1$
and $\omega_2$ are equal when considered as tesselations without
{\it DOE}, i.e. $\omega_1 (\tau_{*} )=\omega_2 (\tau_{*} )$. Thus,
$\omega =\omega_2^{-1}\circ\omega_1 =\gamma_2^{-1}\circ\gamma_1$
is a geometric composition of Whitehead homeomorphisms leaving
invariant the tesselation $\tau_{*}$ without {\it DOE}. Let $K$ be
the intersection of all finite index subgroups arising in this
geometric composition of Whitehead moves. By Relation 4), $\omega$
is uniquely equivalent to a geometric composition of Whitehead
homeomorphisms invariant under the fixed finite index subgroup $K$
of $PSL_2(\mathbb{Z} )$.   Arguing as in Theorem B.4 for the
surface ${\mathbb D}/K$, we conclude that $\omega$ may be reduced
using Relations 1-3) to a composition of automorphisms. $\Box$

\vskip .3 cm

We do not know all the
relations among the generators in Theorem 8.5 and  have only the weaker
statement about normal forms in Theorem 8.7.
One difficulty is that  the normal form for a non-geometric composition of Whitehead
homeomorphisms is not known.

\section{Weil-Petersson two form}

For an oriented smooth punctured surface of finite type with fixed ideal
triangulation $\tau$, it was shown \cite[Appendix A]{P3} that the
Weil-Petersson two form on the Teichm\"uller space pulls-back to the following form

\begin{equation}
-2\sum_{T}(d\log a\wedge d\log b +d\log b\wedge d\log c +
d\log c\wedge d\log a),
\end{equation}
on the decorated Teichm\"uller space $\tilde{T}(F)$, where the
sum is over all triangles $T$ complementary to $\tau $ in $F$ and $(a,b,c)$ are
the edges of $T$ in the correct counter-clockwise cyclic order determined by the
orientation of the surface.

We introduce the Weil-Petersson two form on $\tilde{T}(\S_G )$ by
appropriately normalizing the above expression.  Recall that
$\tau_{*}$ is the Farey tesselation on ${\mathbb D}$. We use the same
notation $\tau_{*}$ for the induced tesselation on the baseleaf of
$\S_G= {\mathbb D}\times\hat{G}/G$ and its canonical TLC
extension to $\S_G$. A fundamental polygon $P$ for the once-punctured
torus group $G$ on ${\mathbb D}$ consists of two ideal triangles $T_1=(e_0
,e_{1/2} ,e_{2/1})$ and $T_2=(e_0 ,e_{-2/1} ,e_{-1/2})$.

We define the tangent vectors on $\tilde{T}(\S_G )$ in terms of
its lambda length parametrization. Note that $Cont^G(\hat{G}
,{\mathbb R}_{>0}^{\tau_{*}})$ is a subset of vector space $Cont^G(\hat{G}
,{\mathbb R}^{\tau_{*}})$ of all continuous functions from $\hat{G}$ onto
${\mathbb R}^{\tau_{*}}$ that are invariant under the action of
$G$. The norm on $Cont^G(\hat{G} ,{\mathbb R}^{\tau_{*}})$ is given by
$$
\| v\| :=\sup_{t\in\hat{G},\ e\in\tau_{*}} |v(t,e)|,
$$
for $v\in Cont^G(\hat{G} ,{\mathbb R}^{\tau_{*}})$. With respect to
this norm, $Cont^G(\hat{G} ,{\mathbb R}_{>0}^{\tau_{*}})$ is an open
subset of the vector space $Cont^G(\hat{G} ,{\mathbb R}^{\tau_{*}})$.
The subspace topology on $Cont^G(\hat{G} ,{\mathbb R}_{>0}^{\tau_{*}})$
induced by the norm on $Cont^G(\hat{G} ,{\mathbb R}^{\tau_{*}})$ coincides
with the compact-open topology. Therefore, the tangent
space at any point of
$\tilde{T}(\S_G )$ is identified with $Cont^G(\hat{G} ,{\mathbb
R}^{\tau_{*}})$.

We define the Weil-Petersson two form on the tangent space at an
arbitrary point $\lambda\in Cont^G(\hat{G} ,{\mathbb
R}_{>0}^{\tau_{*}})$. Let $u,v\in Cont^G(\hat{G} ,{\mathbb
R}^{\tau_{*}})$ be two arbitrary vectors at the tangent space at
$\lambda$. Furthermore, let $\eta^i_{\lambda ,u,v}(t)$ be the
evaluation of the two form $-2(d\log a\wedge d\log b +d\log
b\wedge d\log c + d\log c\wedge d\log a)$ above on the triangle
$T_i\times\{ t\}$ and vectors $u(t,\cdot )$ and $v(t,\cdot )$, for
$i=1,2$ and $t\in\hat{G}$. We define the Weil-Petersson two
form by

$$
\omega_{\lambda} (u,v)=\int_{\hat{G}} [ \eta^{1}_{\lambda
,u,v}(t)+\eta^2_{\lambda ,u,v}(t)] dm(t),
$$
where $m$ is the Haar measure on $\hat{G}$. The actual formulas
are $$ \eta^{1}_{\lambda ,u,v}(t)=-2(\frac{u(t,e_0)
v(t,e_{1/2})}{\lambda (t,e_0)\lambda (t,e_{1/2})}
+\frac{u(t,e_{1/2}) v(t,e_{2/1})}{\lambda (t,e_{1/2})\lambda
(t,e_{2/1})} +\frac{u(t,e_{2/1}) v(t,e_{0})}{\lambda
(t,e_{2/1})\lambda (t,e_0)}),$$ $$\eta^2_{\lambda
,u,v}(t)=-2(\frac{u(t,e_0) v(t,e_{-2/1})}{\lambda (t,e_0)\lambda
(t,e_{-2/1})} +\frac{u(t,e_{-2/1}) v(t,e_{-1/2})}{\lambda
(t,e_{-2/1})\lambda (t,e_{-1/2})} +\frac{u(t,e_{-1/2})
v(t,e_{0})}{\lambda (t,e_{-1/2})\lambda (t,e_0)}).$$ The integral
is over $T_i\times\hat{G}$, for $i=1,2$. In the case when
$\lambda$ is a lift of a decoration on the punctured torus ${\mathbb
D}/G$ and the tangent vectors represent deformations in
$\tilde{T}({\mathbb D}/G )$, the two form $\omega_{\lambda}(u,v)$ is
equal to the classical Weil-Petersson two form as in \cite[Theorem
A.2]{P3} (cf. Theorem~B.2). If
$\lambda$ comes from a lift on some higher genus surface and the
tangent vectors represent TLC deformations of the decorations in
the decorated Teichm\"uller space of that surface, then the two
form is a certain positive multiple of the classical
Weil-Petersson form.

Since $\lambda (t,e)$, for fixed $e\in\tau_{*}$, is a continuous
positive function in $t\in\hat{G}$, it is bounded below away from
zero. Moreover, $u(t,e)$ and $v(t,e)$ are continuous in
$t\in\hat{G}$, for fixed $e\in\tau_{*}$. Therefore,
$\eta^i_{\lambda ,u,v}(t)$ are bounded and continuous real
functions in $t$. Thus, the Weil-Petersson two form
$\omega_{\lambda}(u,v)$ is well-defined, namely, the integral
converges.

\vskip .3 cm

We restrict attention to the base leaf. Let $\eta_{\lambda
,u,v}|_T $ be the evaluation of the two form $-2(d\log a\wedge
d\log b +d\log b\wedge d\log c + d\log c\wedge d\log a)$ on the
triangle $T=(a,b,c)$ in ${\mathbb D} -\tau_{*}$, where $d\log a
=\dot{a}/a$, etc... are given their corresponding values in terms
of $\lambda (id,a)$, $u(id,a)$, and so on,  on the baseleaf
${\mathbb D}$. We consider the average evaluation of $\eta$ on the
baseleaf. Consider a sequence of characteristic  subgroups $G_n$ ( cf. Section 2) of
$PSL_2\mathbb{Z}$ such that $\cap_{n} G_n=\{ id\}$. Let $P_n$ be a
sequence of fundamental polygons for $G_n$ such that the frontier edges of
$P_n$ belong to the Farey tesselation $\tau_{*}$, and the $P_n$ are
nested.

\vskip .3 cm

\paragraph{\bf Proposition 9.1}{\it
Let $\lambda\in Cont^G(\hat{G} ,{\mathbb R}_{>0}^{\tau_{*}} )$ and let
$u,v\in Cont^G(\hat{G} ,{\mathbb R}^{\tau_{*}} )$ be two tangent
vectors based at $\lambda$.  Then

$$
\omega_{\lambda} (u,v)=\lim_{n\to\infty} \frac{2}{k(n)}\sum_{T\in
P_n} \eta_{\lambda ,u,v} |_T,
$$
where the sum is over all triangles $T$ complementary to
$\tau_{*}$ in the fundamental polygon $P_n$, and $k(n)$ is the
number of triangles in $P_n$.}

\vskip .3cm

\paragraph{\bf Proof} We first show that $\lim_{n\to\infty}
\frac{2}{k(n)}\sum_{T\in P_n}\eta_{\lambda ,u,v} |_T$ does not
depend on the choice of the fundamental polygons $P_n$. Indeed,
for $n$ large enough the quantities $\lambda$, $u$ and $v$ are
almost invariant under the group $G_n$. Thus, the quantity
$\eta_{\lambda ,u,v}|_T$ is also almost invariant under $G_n$.
Namely, given $\epsilon >0$ there exists $n$ such that
$|\eta_{\lambda ,u,v}|_T -\eta_{\lambda ,u,v}|_{\gamma
T}|<\epsilon$ for all triangles $T$ complementary to $\tau_{*}$
whenever $\gamma\in G_n$. Hence, if the limit exists it does not
depend on the choice of $P_n$.

For a fixed sequence $P_n$ of fundamental polygons, the sum
$\frac{2}{k(n)}\sum_{T\in P_n} \eta_{\lambda ,u,v} |_T$ is a
Cauchy sequence. This immediately follows from the above
observation of almost invariance of $\eta_{\lambda ,u,v} |_T$.

It remains to show the equality between $\omega_{\lambda} (u,v)$
and the limit of the sum. In the case when $\lambda$, $u$ and $v$
are TLC, the equality follows from the definition of the Haar
measure and an observation that our normalization factor
$\frac{2}{k(n)}$ is chosen correctly. The general case follows from
the uniform continuity of $\eta^i_{\lambda ,u,v}(t)$ and the fact
that $\eta_{\lambda ,u,v} |_{\gamma T} =\eta^i_{\lambda ,u,v}(t)$
in the case when $t=\gamma\in\hat{G}$ is an element of $G<\hat{G}$
and $T$ is the corresponding triangle $T_i\times\{ \gamma\}$.
$\Box$

\vskip .3 cm

We show that the Weil-Petersson
two form pushes forward from $\tilde{T} (\S_G )$ onto the Teichm\"uller
space
$T(\S_G )$. Consider the projection map $\pi :\tilde{T} (\S_G )\to
T(\S_G )$. To show that there exists a well-defined push-forward
of the two form, it is enough to show that $\omega_{\lambda^1}
(u^1,v^1)=\omega_{\lambda^2}(u^2,v^2)$ whenever $\pi (\lambda^1
)=\pi (\lambda^2)$, $d\pi_{\lambda^1}(u^1)=d\pi_{\lambda^2}(u^2)$
and $d\pi_{\lambda^1}(v^1)=d\pi_{\lambda^2}(v^2)$.

\vskip .3 cm

\paragraph{\bf Theorem 9.2}
{\it The Weil-Petersson two form on $\tilde{T}(\S_G )$ projects to
a non-degenerate two form on the Teichm\"uller space $T(\S_G )$.}

\vskip .3cm

\paragraph{\bf Proof}
Suppose that $\lambda_i$, $u_i$ and $v_i$, $i=1,2$ are as above.
If they are invariant under a finite index subgroup $K$ of $G$,
then $\omega_{\lambda_1} (u_1,v_1)=\omega_{\lambda_2} (u_2,v_2)$.
Indeed, since the Weil-Petersson form in this case is a multiple
of the Weil-Petersson form on $\tilde{T}({\mathbb D}/K)$ (which is the
lift of the Weil-Petersson form on $T({\mathbb D}/K)$ by Theorem~B.2), it
projects to the multiple of the Weil-Petersson form on $T({\mathbb
D}/K)$. Thus, it satisfies $\omega_{\lambda_1}
(u_1,v_1)=\omega_{\lambda_2} (u_2,v_2)$. In general, it is enough to show
that any $\lambda_i$, $u_i$ and $v_i$,
$i=1,2$ can be approximated by TLC objects which satisfy the property in
the paragraph above the theorem. Indeed, the union of decorated
Teichm\"uller spaces of all finite surfaces is dense in the
decorated Teichm\"uller space $\tilde{T} (\S_G )$ and the
inclusion is smooth. It remains to see that two tangent vectors
$u_1,u_2$ which project to the same vector in the tangent space of
$T(\S_G )$ can be approximated by two TLC vectors $u_1^c ,u_2^c$
which project to the same TLC vector in the tangent space of
$T(\S_G )$. We postpone this to Lemma 9.3. Since the formula for
$\omega_{\lambda} (u,v)$ is continuous in its entries, the limit
also satisfies the desired push-forward property. Thus, the
Weil-Petersson two form indeed projects to the Teichm\"uller space
$T(\S_G )$.

It remains to show that the push-forward of the Weil-Petersson two
form is non-degenerate. To see this, we must show that for
an arbitrary vector $v$ at a point $\lambda\in\tilde{T} (\S_G )$
which projects to a non-zero vector on $T(\S_G )$ there exists
another vector $u$ such that $\omega_{\lambda} (u,v)\neq 0$. Note
that a vector $v\in Cont^G(\hat{G},\mathbb{R}^{\tau_{*}})$
projects to a trivial vector on $T(\S_G )$ if and only if it is
tangent to a path in $\lambda_s\in
Cont^G(\hat{G},\mathbb{R}^{\tau_{*}}_{>0})$ such that
$cr(\lambda_s (t,a),\lambda_s
(t,b),\lambda_s(t,c),\lambda_s(t,d))=const$ for all quadrilaterals
$(a,b,c,d)$ which are union of two adjacent triangles of the Farey
tesselation $\tau_{*}$ and for all $t\in\hat{G}$, where $cr$ denotes the
cross ratio. Thus,
$v$ projects to a trivial vector on $T(\S_G )$ if and only
if there exists a rectangle $Q=(a,b,c,d)$ consisting of two
adjacent triangles in ${\mathbb D}-\tau_{*}$ such that
$|v(t,a)/\lambda (t,a)+v(t,c)/\lambda (t,c)-v(t,b)/\lambda
(t,b)-v(t,d)/\lambda (t,d)|\geq \epsilon
>0$ for $t\in\hat{K}$, where $K$ is a finite index subgroup of
$G$. Thus, either $|v(t,a)/\lambda (t,a)+v(t,c)/\lambda
(t,c)|\geq\epsilon /2$ or $|v(t,b)/\lambda (t,b)-v(t,d)/\lambda
(t,d)|\geq \epsilon /2$, for $t\in\hat{K}_1$, where $K_1$ is a
finite index subgroup of $K$.  We may assume that they have the
same sign for $K_1$ small enough and that $|v(t,a)/\lambda
(t,a)+v(t,c)/\lambda (t,c)|\geq\epsilon /2$, for $t\in\hat{K}_1$. Let
$e\in\tau_{*}$ be the diagonal of $Q$, and choose $u$ such that
$u(t,e)\equiv 1$ for
$t\in\hat{K}_1$ and equals zero otherwise.

Since
$\omega_{\lambda}(u,v)=-2\int_{ \hat{K}_1}\frac{1}{\lambda (t,e)}
(\frac{v(t,a)}{\lambda (t,a)} +\frac{v(t,c)}{\lambda
(t,c)})dm(t)$ and either $\frac{1}{\lambda (t,e)}
(\frac{v(t,a)}{\lambda (t,a)} +\frac{v(t,c)}{\lambda
(t,c)})\geq\epsilon /2$ or $\frac{1}{\lambda (t,e)}
(\frac{v(t,a)}{\lambda (t,a)} +\frac{v(t,c)}{\lambda (t,c)})\leq
-\epsilon /2$, for all $t\in\hat{K}_1$, we conclude that
$\omega_{\lambda} (u,v)\neq 0$ as required. $\Box$

\vskip .3 cm

\paragraph{\bf Lemma 9.3} {\it Let $\epsilon >0$ and let $u_1,u_2$
be two tangent vectors at $\lambda_1,\lambda_2$, respectively. If
$d\pi_{\lambda_1} (u_1)=d\pi_{\lambda_2} (u_2)$ then there exist
 two TLC points $\lambda_1^c,\lambda_2^c$ and two TLC
vectors $u_1^c,u_2^c$ such that $d\pi_{\lambda_1^c} (u_1^c)=d\pi_{\lambda_2^c}
(u_2^c)$, $\|\lambda_i-\lambda_i^c\|_{\infty} <\epsilon$ and $\|
u_i-u_i^c\|_{\infty} <\epsilon$, for $i=1,2$.}

\vskip .3cm

\paragraph{\bf Proof} By our assumption, there exist two paths
$\lambda_{1,s}$ and $\lambda_{2,s}$ in $\tilde{T} (\S_G
)$ such that $\lambda_{1,0}=\lambda_1$, $\lambda_{2,0}=\lambda_2$,
$\pi (\lambda_{1,s})=\pi (\lambda_{2,s})$, $\frac{d}{ds}
\lambda_{1,s} |_{s=0}=u_1$ and $\frac{d}{ds} \lambda_{2,s} |_{s=0}
=u_2$. Let $n$ be an nteger large enough that $\lambda_{1,s}$
and $\lambda_{2,s}$ on the baseleaf are almost invariant under
$G_n$. We form two new paths of TLC lambda lengths which satisfy
the required properties. We recall (cf. Lemma~A1iv) that the
projection $\pi :\tilde{T} (\S_G )\to T(\S_G )$ can be given by
the formula
\begin{equation}
\label{one} cr(Q)=\frac{ac}{bd},
\end{equation}
where $cr(Q)$ is a cross-ratio of the endpoints of a quadrilateral
$Q=T_1\cup T_2$ for two adjacent triangles $T_1,T_2\in {\mathbb
D}-\tau_{*}$ and $a,b,c,d$ are lambda lengths evaluated at consecutive
edges of $Q$. A necessary and sufficient condition that $\pi
(\lambda_{1,s})=\pi (\lambda_{2,s})$ is thus to have
\begin{equation}
\label{two}
\frac{\lambda_{1,s}(a)\lambda_{1,s}(c)}{\lambda_{1,s}(b)\lambda_{1,s}(d)}
=
\frac{\lambda_{2,s}(a)\lambda_{2,s}(c)}{\lambda_{2,s}(b)\lambda_{2,s}(d)}
\end{equation}
for every such $Q=(a,b,c,d)$. We form new paths:
$$
\lambda_{i,s}^c (a)=\lim_{k\to\infty}(\Pi_{\gamma\in H_k\subset
G_n}\lambda_{i,s}(\gamma a))^{\frac{1}{k}},
$$
where $H_k$ contains $k$ elements of $G_n$, $H_k\subset H_{k+1}$
and $\cup_{k=1}^{\infty} H_k =G_n$. The limits exist, they do not
depend on the choice of $H_k$ and are TLC. Moreover,
$\lambda_{i,s}^c$ satisfy (\ref{two}) because at each finite stage,
$(\Pi_{\gamma\in H_k\subset G_n}\lambda_{i,s}(\gamma
a))^{\frac{1}{k}}$ satisfy (\ref{two}). They are close to the
original paths for $n$ large, and their tangent vectors at $s=0$
approximate the original tangent vectors, are also TLC, and
their projections agree. The lemma follows. $\Box$

\vskip .3 cm

\paragraph{\bf Theorem 9.4} {\it The Weil-Petersson two form on
$\tilde{T} (\S_G )$ is invariant under the modular group
$Mod_{BLP} (\S_G )$. Therefore, it descends to $\tilde{T} (\S_G
)/Mod_{BLP} (\S_G )$ and $T(\S_G )/Mod_{BLP} (\S_G )$.}

\vskip .3cm

\paragraph{\bf Proof} We recall from Theorem~8.5 that the modular group
is generated by equivariant Whitehead moves and $PSL_2(\mathbb{Z}
)$. It is therefore enough to show the invariance under these
elements of $Mod_{BLP}(\S_G )$. It is a calculation in \cite{P3}
(that the reader could provide: show that the two form is algebraically
invariant under Ptolemy transformations Lemma A.1(iii)) that a Whitehead move along a
single edge does not change the sum
$\eta_{\lambda ,u,v} |_{T_1} +\eta_{\lambda ,u,v} |_{T_2}$ over
the two triangles which contain the edge on their boundaries.
Since an equivariant Whitehead move is decomposed into Whitehead
moves along disjoint edges, the invariance follows. The action of
$PSL_2(\mathbb{Z} )$ fixes $\tau_{*}$ and changes the marking. It
is obvious that the Weil-Petersson two form is independent of the
marking. The invariance under the whole $Mod_{BLP}(\S_G
)$ follows. $\Box$

\vskip .3cm

\section{Concluding Remarks and Questions}

\vskip .2cm

To what extent is there number theory embedded in our
constructions? The ``center'' of ${\mathcal C}(\tau )$, i.e., a
point corresponding to equating all lambda lengths, covers an
arithmetic punctured surface for each TLC tesselation $\tau$
\cite[Section 6]{P1}, which is evidently related to Grothendieck's
dessins d'enfant \cite{LS}.  In fact, equivariant Whitehead moves were considered
before in this context  \cite{P4}.  Is there any connection between
the absolute Galois group and the full (non-BLP) mapping class group
of the solenoid? Furthermore, just as the Euclidean
solenoid \cite{Od} is related to the rational adeles, our
constructions give an adelic type structure to the punctured
solenoid itself. How might the Teichm\"uller theory developed
here, a kind of deformation  theory of arithmetic punctured
surfaces be manifest algebraically? Furthermore, the calculation
of the index of the characterstic subgroups $G_n$ of
$G=PSL_2({\mathbb Z})$ is an open problem in number theory; on one
hand, it is evidently related to fatgraph enumeration
\cite[Appendix B]{P3}, and on the other, the results of Section 9
suggest connections with Weil-Petersson volumes.

More generally, the regularization of Weil-Petersson form in Section 9
together with the strong topology on the lambda length functions seems to
give a much more satisfactory universal geometry than
\cite{P1}.  With this symplectic or Poisson structure, the
global lambda length coordinates, and the generators of the modular
group for the solenoid, all the ingredients are in place for a Kashaev \cite{Ka} or
Chekhov-Fock type quantization \cite{CF} of the decorated Teichm\"uller
theory developed here.

Are $\stackrel{\circ}{\mathcal C}(\tau )$ or ${\mathcal C}(\tau )$
contractible for a TLC tesselation or paving $\tau$?  The classical
tools (from \cite{P2} used to prove Lemma~B.4) seem to be unavailable
here.  Which non-TLC pavings  of ${\mathbb D}$ arise
from the convex hull construction?

We give only an infinite set of generators, and we ask if
$Mod_{BLP}(\S)$ is finitely generated.
Of course, we would also hope to give a presentation of this group,
beyond the normal forms presented here; perhaps one can mimic the
case of punctured surfaces by first proving contractibility of the ${\mathcal C}(\tau )$.
It also seems possible to
characterize the mapping class like elements of $Mod_{BLP}(\S)$,
i.e., those homeomorphisms of $\S$ that arise from lifts of
homeomorphisms of punctured surfaces,  in terms of our generators,
and we wonder if this might be useful to address a question from
\cite{Od}: do the mapping class like homeomorphisms generate the
modular group?

\vskip .3 cm

\section*{Appendix A-Pinched lambda lengths and the convex hull
construction}

\vskip .2in

Let ${\mathbb D}$ denote the unit disk with boundary $S^1$.  A
{\it tesselation} $\tau$ of ${\mathbb D}$ is a locally finite
decomposition of ${\mathbb D}$ into ideal triangles.  We shall
think of $\tau$ as a set of edges and let $\tau ^\infty\subseteq
S^1$ denote the collection of ideal vertices of $\tau$.  In
particular, for the Farey tesselation $\tau _*$, we have $\tau
_*^\infty={\mathbb Q}\cup\{\infty\}=\bar{\mathbb Q}\subseteq S^1$.
It will be useful sometimes to specify a distinguished oriented
edge or {\it DOE} of a tesselation. In particular, the standard
{\it DOE} of the Farey tesselation is the edge connecting $-1$ to
$+1$.

Define Minkowski three-space to be the vector space ${{\mathbb R}}^3$
with indefinite pairing $<\cdot ,\cdot >$ whose quadratic form is
$x^2+y^2-z^2$ in the usual coordinates.  The upper-sheet of the
hyperboloid is ${\mathbb H}=\{w=(x,y,z):
<w,w>=-1,~z>0\}$, and the positive light cone in
Minkowski space is
$L^+=\{ w=(x,y,z):<w,w>=0,~z>0\}$.  The former is a model for the
hyperbolic plane, where the distance $\Delta$ between two points
$u,v\in{\mathbb H}$ is given by
$\Delta^2=~{\rm cosh}^2 <u,v>$.  The latter
is identified with the space of all horocycles via the
correspondence
$L^+\ni w\mapsto
\{ u\in{\mathbb H}: <w,u>=-1\}\subseteq {\mathbb H}$, and there is the
corresponding identification $L^+/{{\mathbb R}}_{>0}\equiv S^1$ which
maps a horocycle onto its center.  The {\it lambda length} of
$u,v\in L^+$ is defined to be
$\Lambda (u,v)=\sqrt {-<u,v>}$, and geometrically it is $\sqrt{2~{\rm
exp}~\delta}$, where $\delta$ is the signed hyperbolic distance between
the corresponding horocycles, taken with positive sign if and only if the
horocycles are disjoint.  An affine plane in Minkowski space is
respectively elliptic, parabolic, hyperbolic if and only if it determines
the corresponding conic section, or equivalently if and only if the
Minkowski normal $v$ to the plane $\{ w: <w,v>=-1\}$ lies interior
to, lies on, or lies exterior to the light-cone.

\vskip .2in

\noindent {\bf Lemma A.1}~~\it Given pairwise non collinear points
$u_1,u_2,u_3,u_4\in L^+$ so that the rays determined by $u_1,u_3$
separate those determined by $u_2,u_4$, set $\lambda
_{ij}=\sqrt{-<u_i,u_j>}$.

\vskip .1in

\noindent {\rm (i)}~\cite[Lemmas 2.3, 2.4]{P2} ~~Given three
positive real numbers
$\ell _{12},\ell_{13},\ell_{23}\in{{\mathbb R}}_{>0}$, there are unique
points $v_i$ in the ray in $L^+$ determined by $u_i$, for
$i=1,2,3$, so that $\ell _{ij}=\lambda _{ij}$.  Furthermore, given
$v_1,v_2\in L^+$ so that $\ell _{12}=\lambda _{12}$, there is a
unique $v_3\in L^+$ on either side of the plane through the origin
determined by $v_1,v_2$ so that also $\ell _{13}=\lambda _{13}$
and $\ell _{23}=\lambda _{23}$. In each case, the points $v_i$
depend continuously on the lambda lengths.

\vskip .1in

\noindent{\rm (ii)}~\cite[Lemma 2.2]{P2}~~The plane
spanned by
$u_1,u_2,u_3$ is elliptic if and only if all three strict triangle
inequalities hold amongst
$\lambda _{12},\lambda _{13},\lambda _{31}$, and it is parabolic
if and only if a triangle equality holds.

\vskip .1in

\noindent {\rm (iii)}~\cite[Proposition~2.6a]{P2}~~
The {\rm Ptolemy equation} holds:
$$\lambda _{13}\lambda_{24}= \lambda _{12}\lambda_{34}+\lambda
_{14}\lambda_{23}.$$

\vskip .1in

\noindent {\rm (iv)}~\cite[Lemma A.2]{P3}~~Let $\bar
u_i\in S^1$ denote the projection of
$u_i\in  L^+$, for $i=1,2,3,4$.  The cross-ratio is given as follows:
normalizing so that
$\bar u_1\mapsto 1$,
$\bar u_2\mapsto 0$, $\bar u_4\mapsto\infty$, then $\bar u_3\mapsto
{{\lambda_{23}\lambda_{34}}\over{\lambda_{12}\lambda_{14}}}.$

\vskip .1in

\noindent {\rm (v)}~~\cite[Proposition 2.6b]{P2}~~The signed volume of
the Euclidean tetrahedron spanned by
$u_1,u_2,u_3,u_4$ is given by
$(\lambda _{12}\lambda_{23}\lambda _{34}\lambda_{41})~\sigma$, where the
{\rm simplicial coordinate} $\sigma$ is defined by:
$$\sigma = {{\lambda _{12}^2 +\lambda_{23}^2-\lambda _{31}^2}\over{
\lambda _{12}\lambda_{23}\lambda _{31}}} +{{\lambda _{14}^2
+\lambda_{43}^2-\lambda _{31}^2}\over{ \lambda
_{14}\lambda_{43}\lambda _{31}}},$$ and the sign is positive if
and only if the edge connecting $u_1,u_3$ lies below the edge
connecting $u_2,u_4$.\rm

\vskip .2in

We claim that any function $\Lambda :\tau _*\to{{\mathbb R}}_{>0}$
gives rise to a unique $h:\tau _*^\infty\to L^+$ realizing the
putative lambda lengths in the sense that for all $e\in\tau _*$,
we have $\Lambda (e)=\sqrt{-<h(u),h(v)>}$, where $u,v$ are the
endpoints of $e$.  To see this, we may uniquely realize the
putative lambda lengths on the triangle to the left of the {\it
DOE} by points in the rays in $L^+$ lying over $\pm 1,-i$ by the
first statement in part (i) of Lemma~A.1.  We may then use the
second statement in part (i) to recursively define the required
function $h:\tau _*^\infty\to L^+$, where at each stage in the
construction in the notation of part (i), $v_3$ lies on the other
side of the plane through the origin containing $v_1,v_2$ from the
triangle to the left of the {\it DOE}.  Post-composing with the
projection, we define the ``characteristic map'' $\bar h:\tau
_*^\infty\to L^+/{{\mathbb R}}_{>0}\equiv S^1$ and can ask whether
$\bar h$ interpolates a homeomorphism $S^1\to S^1$, and if so,
what is the nature of this homeomorphism.

We say that $\Lambda :\tau _*\to{{\mathbb R}}_{>0}$ is {\it pinched} if
there is a constant $M>1$ so that for all $e\in\tau _*$, we have
$M^{-1}<\Lambda (e)<M$.

\vskip .2in

\noindent {\bf Theorem A.2}~~\cite[Theorems 6.3 and 6.4]{P1}~~
\it ~~If
$\Lambda :
\tau _*\to {{\mathbb R}}_{>0}$ is pinched, then $h(\tau _*^\infty
)\subseteq L^+$ satisfies the following two properties:

\vskip .1in

\leftskip .3in

\noindent {\rm (1)}~$h(\tau _*^\infty )$ is {\rm discrete}, i.e.,
below any elliptic plane in Minkowski space, there are only
finitely many points of $h(\tau _*^\infty )$.

\vskip .1in

\noindent {\rm (2)}~$h(\tau _*^\infty )$ is {\rm radially dense},
i.e., in the union of any open set of rays in $L^+$, there is some
point of $h(\tau _*^\infty )$.

\vskip .1in

\leftskip=0ex

\noindent Furthermore, $\bar h:\tau _*^\infty\to S^1$ interpolates
a homeomorphism $\phi :S^1\to S^1$ which is quasi-symmetric.\rm

\vskip .2in

If $\Lambda\in{{\mathbb R}}_{>0}^{\tau _*}$ is pinched, let ${\mathcal
B}=h(\tau _*^\infty ) \subseteq L^+$, and let $C\subseteq {{\mathbb
R}}^3$ denote the closed convex hull of ${\mathcal B}$ (in the vector
space underlying Minkowski space).  We may think of the boundary
$\partial C$ of $C$ in Minkowski space as a piecewise-linear
approximation to the upper sheet of the unit hyperboloid with its
vertices in $L^+$:

\vskip .2in

\noindent {\bf Theorem A.3} \cite[Lemma 7.2]{P1} ~~\it ~Suppose
that
${\mathcal B}\subseteq L^+$ is a discrete and radially dense subset with
closed convex hull $C$. Then $C\cap L^+$ is the set of points of
the form $tz$, where $t\geq 1$ and $z\in {\mathcal B}$; each ray
from the origin lying inside $L^+$ meets $\partial C$ exactly
once. The boundary $\partial C$ is the union of $C\cap L^+$ and a
countable collection of codimension one {\rm faces}
$F_1,F_2,\ldots$.  Each such face is the convex hull of some
coplanar subset $X\subseteq{\mathcal B}$.  The affine plane
containing $X$ is either parabolic or elliptic, and if $X$ is
infinite, then this affine plane is parabolic. The set of faces is
locally finite in the interior of $L^+$.\rm

\vskip .2in

An edge of $\partial C$ determines a geodesic in ${\mathbb D}$ in the
natural way, and we let $\tau\subseteq{\mathbb D}$ denote the set of
geodesics corresponding to all the edges of $\partial C$.
According to the previous theorem, $\tau$ is a locally finite
collection of disjoint geodesic whose complementary regions are
either finite-sided polygons (corresponding to elliptic or
parabolic support planes) or infinte-sided (corresponding to
parabolic support planes), and we call such a decomposition of
${\mathbb D}$ a {\it paving}.
Generically, no four points of
${\mathcal B}$ are coplanar, and the paving $\tau$ is a tesselation of
${\mathbb D}$.

We shall adopt this as standard notation: if $\Lambda\in{{\mathbb
R}}_{>0}^{\tau _*}$ is pinched, then the corresponding $h_\Lambda
:\tau _*^\infty\to L^+$ has discrete radially dense (drd) image
${\mathcal B}_\Lambda =h_\Lambda (\tau _*^\infty)\subseteq L^+$
and the projection $\bar h_\Lambda:\tau _*^\infty\to S^1$
interpolates a quasi-symmetric homeomorphism $\phi _\Lambda
:S^1\to S^1$; we may also sometimes extend $\phi _\Lambda$ to a
quasi-conformal map $\phi _\Lambda:{\mathbb D}\to {\mathbb D}$ say using
\cite{DE}. The edges in the boundary of the closed convex hull $C_\Lambda$
of
${\mathcal B}_\Lambda$ project to a locally finite collection
$\tau _\Lambda$ of disjoint geodesics in ${\mathbb D}$.

By the geometric interpretation of simplicial coordinates in
Lemma~A.1(v), it follows that the simplicial coordinate of (the edge of
$C_\Lambda$ corresponding to) an edge of $\tau _\Lambda$ is a
well-defined non-negative real number.  In particular in the generic
case that
$\tau$ is a triangulation, all the simplicial coordinates are
strictly positive.

Given a tesselation $\tau$ and given any unoriented edge $f$ of
$\tau$, define the {\it Whitehead move along} $f$ to be the
tesselation $\tau _f=\tau\cup\{ g\} -\{ f\}$, where $f,g$ are the
diagonals of an ideal quadrilateal with frontier in $\tau$.
Furthermore, if $\tau$ comes equipped with a {\it DOE} $e$, then
there is an corresponding {\it DOE} on $\tau _f$, where if $e=f$,
then the orientations on $e,f$ in this order is consistent with
the orientation of ${\mathbb D}$ itself, and if $e\neq f$, then
the {\it DOE} is unchanged.

Thus, sequences of Whitehead moves act on tesselations or on
tesselations with {\it DOE}.  Sequences of Whitehead moves on
interior edges also act on triangulations of finite polygons and
on triangulations with {\it DOE} of finite polygons.

\vskip .2in

\noindent{\bf Theorem A.4}~\cite[Lemma 4.4]{P1}~~\it Fix a
finite-sided connected polygon and consider two triangulations of
it with {\it DOE}.  Then there is a finite sequence of Whitehead
moves from one to the other, i.e., sequences of Whitehead moves
act transitively on tesselations with {\it DOE} of a finite
polygon.\rm

\vskip .3 cm

\section*{Appendix B-Punctured surfaces}

\vskip .2in

This appendix is intended to recall the principal results
regarding punctured surfaces from \cite{P2} and to give an
explicit presentation for the modular groups of punctured
surfaces.  Fix a surface $F$ with negative Euler characteristic
and $s\geq 1$ many punctures.  An ``ideal triangulation'' of $F$
is the homotopy class of a family $\tau$ of disjointly embedded
arcs connecting punctures so that each component of $F-\cup\tau$
is a collection of triangles.  A ``decoration'' on $F$ is the
specification of one horocycle centered at each puncture of $F$.
The forgetful map $\tilde{T}(F)\to\T(F)$ from decorated to
undecorated Teichm\"uller space of $F$ is a principal ${\mathbb
R}_{>0}^s$-bundle, where we may take for instance the hyperbolic
lengths of the horocycles as coordinates on ${\mathbb R}_{>0}^s$.
As a point of notation, we shall let $\Gamma\in T(F)$ denote the
(class of a) Fuchsian group uniformizing a point of $T(F)$, and
$\tilde\Gamma\in\tilde T(F)$ a decorated hyperbolic structure with
underlying $\Gamma\in T(F)$.

\vskip .2cm

\noindent {\bf Theorem B.1}~~\cite[Theorem 3.1]{P2} ~~\it Fix an
ideal triangulation $\tau$ of $F$.  Then the natural map $\tilde
T(F)\to {\mathbb R}_{>0}^\tau$, which associates to
$\tilde\Gamma\in\tilde T(F)$ the function that assigns the lambda
length of $e\in\tau$ for $\tilde\Gamma$, is a homeomorphism onto.
That is, lambda lengths on the edges of $\tau$ give a
parametrization of $\tilde T(F)$.\rm

\vskip .2cm

The proof was sketched after Lemma~A.1: use lambda lengths on the
lift of $\tau$ to the universal cover to recursively define
$h:\tau _*^\infty\to L^+$.  This characteristic map can easily be
shown directly to interpolate a homemorphism of the circle (or
alternatively use Theorem~A.1 since the finitely many lambda
lengths are {\it a priori} pinched).  Thus, there is a
corresponding tesselation of ${\mathbb D}$, which is invariant
under a Fuchsian group by Poincar\'e's fundamental polygon
theorem.  This gives the point $\Gamma\in T(F)$, and the
construction furthermore gives a decoration $\tilde\Gamma\in\tilde
T(F)$ on it.   This provides the inverse homeomorphism to the map
in Theorem~B.1.

\vskip .2cm

\noindent {\bf Theorem B.2}~~\cite[Theorem A.1]{P3} ~~\it In the
notation of Theorem~B.1, the Weil-Petersson K\"ahler two from on
$T(F)$ lifts to
$$-2\sum ~{\rm d}ln a~\wedge ~{\rm d}ln b
+~{\rm d}ln b~\wedge ~{\rm d}ln c+~{\rm d}ln c~\wedge ~{\rm d}ln
a$$ on $\tilde T(F)$, where the sum is over all triangles
complementary to $\tau$ in $F$ with oriented boundary having lamba
lengths $a,b,c$ in this counter-clockwise cyclic order.\rm

\vskip .2cm

The proof is a calculation in \cite[Appendix A]{P3} starting from
Wolpert's formula.

\vskip .2cm

Turning now to the convex hull construction in this case, the drd
set ${\mathcal B}$ consists of $s\geq 1$ many $\Gamma$-orbits in
the light cone $L^+$.  The support planes of the convex hull $C$
of ${\mathcal B}$ are all elliptic (cf. \cite[Proposition
4.4]{P2}), hence finite sided by discreteness.  This is in
contrast to Theorem~A.2, where parabolic support planes can occur.
The extreme edges of the $\Gamma$-invariant convex body $C$
project to ${\mathbb D}$ as usual and then to a family $ \tau _
{\tilde\Gamma }$ of arcs connecting punctures in $F$.

\vskip .2cm

\noindent{\bf Easy Lemma B.3}~~\cite[Theorem 4.5]{P2}~~\it For
every $\tilde\Gamma\in\tilde T(F)$, $\tau _ {\tilde\Gamma }$ is a
paving, i.e., consists of disjointly embedded arcs, no two of
which are homotopic, connecting punctures so that each
complementary region is a polygon.\rm

\vskip .2cm

For any paving $\tau$ of $F$, define
\begin{equation*}
\begin{array}l
\stackrel{\circ}{\mathcal C}(\tau )=\{\tilde\Gamma \in \tilde
T({F}): \tau_{\tilde\Gamma}=\tau\} \\ \ \ \cap\\{\mathcal C}(\tau
)=\{\tilde\Gamma \in \tilde T({F}): \tau _{
\tilde\Gamma}\subseteq\tau\}\\ \ \ \cap\\ \tilde T(F).
\end{array}
\end{equation*}

\vskip .2cm

\noindent{\bf Hard Lemma B.4}~~\cite[Theorem 5.4]{P2}~~\it Each
${\mathcal C}(\tau )$ has the natural structure of an open simplex
$\stackrel{\circ}{\mathcal C}(\tau )$ plus certain of its
faces.\rm

\vskip .2cm

\noindent Together, the two previous lemmas give:

\vskip .2cm

\noindent {\bf Theorem B.5}~~\cite[Theorem 5.5]{P2} ~~\it $\{
\stackrel{\circ}{\mathcal C}(\tau ):~\tau~{\rm
is~a~paving~of}~F\}$ provides a cell decomposition of $\tilde
T(F)$ which is invariant under the modular group $Mod(F)$ of
$F$.\rm

\vskip .2cm

\noindent {\bf Corollary~B.6}~\cite[Proposition 7.1]{P2}
\it Compositions of Whitehead moves act transitively on the
collection of all ideal triangulations of a fixed punctured
surface.\rm

\vskip .2cm

\noindent {\bf Proof}~~$\stackrel{\circ}{\mathcal C}(\tau
)\neq\emptyset$ for any ideal triangulation $\tau$ of $F$; to see
this, take all lambda lengths to be unity and compute simplicial
coordinates (cf. Lemma A.1(v)). Since $\tilde T(F)$ is path
connected, there is a path between the top-dimensional cells
corresponding to any two ideal triangulations of $F$.  Putting
this path in general position with respect to the codimension-one
faces of the cell decomposition shows that compositions of
Whitehead moves act transitively on ideal triangulations of $F$.
~$\Box$

\vskip .2cm

Consider the {branching locus} $\{\tilde\Gamma\in\tilde
T(F):\exists\phi\in Mod(F)~{\rm with}~\phi
~\tilde\Gamma=\tilde\Gamma\}$ with projection $\tilde{\mathcal
L}\subseteq \tilde M=\tilde M(F)=\tilde T(F)/Mod(F)$. Recall (cf.
\cite{Thbook},\cite{Igusa}) that the ``orbifold fundamental
group'' $Mod(F)=\pi _1^{orb}(\tilde M)$ is the semi-direct product
of the usual topological fundamental group $\pi _1^{top}(\tilde
M-\tilde{\mathcal L})$ with a finite group determined by the
branching, i.e., by the symmetry groups of the underlying pavings.

\vskip .2cm

\noindent{\bf Corollary B.7}~~\it The modular group $Mod(F)=\pi
_1^{orb}(\tilde M(F))$ of any punctured surface $F$ of negative
Euler characteristic is the stabilizer of any object in the
groupoid with finite presentation:

\vskip .2cm

\noindent generators are provided by  Whitehead moves between
ideal triangulations and the union of automorphism groups of all
ideal triangulations of $F$;

 \vskip .2cm

\noindent relations are given by involutivity, the pentagon and
commutativity  relations (as in Section 8),   and { invariance
under $\phi\in Mod(F)$}, namely, for any ideal triangulation $\tau
$ of $F$, we have:  {\rm i) } the Whitehead move on $e$ in $\tau$
is identified with the Whitehead move on $\phi (e)$ in $\phi (\tau
)$, and likewise, $Aut(\tau )$ is identified with $Aut(\phi (\tau
))$; and {\rm ii)} if $w$ is a sequence of Whitehead moves
beginning at $\tau$ and ending at $\phi (\tau )$,  then up to pre-
and post-composition with elements of $Aut(\tau )\approx Aut(\phi
(\tau ))$, we identify $\phi$ with $w$.\rm

\vskip .2cm

\noindent {\bf Proof}~~That the putative relations hold amongst
the putative generators follows from the discussion in Section 8
plus the obvious invariance i) and ii) under $Mod(F)$.  As in
\cite{Igusa}, Whitehead moves together with the automorphism
groups of all pavings of $F$ generates $Mod(F)$ since compositions
of Whitehead moves act transitively on ideal triangulations of a
polygon (by Theorem A.4). Furthermore, if $\phi\in Aut (\sigma )$
for some paving $\sigma$, then we may arbitrarily extend to an
ideal triangulation $\tau \supset\sigma$.  By Corollary~B.6, there
is some sequence $w$ of Whitehead moves from $\phi (\tau )$ to
$\tau$, and by naturality condition ii), the composition given by
first $\phi$ and then $w$ must lie in $Aut(\tau )$.  It follows
that Whitehead moves and automorphism groups of ideal
triangulations alone indeed generate $Mod(F)$.

As to relations, a homotopy of loops or paths in the orbifold
$\tilde M$ gives rise to a homotopy also in the underlying space,
and this homotopy may be put into general position with the
codimension two skeleton of the cell decomposition.  The orbifold
fundamental group is a semi-direct product of $\pi _1^{top}(\tilde
M-\tilde{\mathcal L})$ with a finite group as above, and this
finite group is contained in the span of the generators by the
previous paragraph.  Up to the action of this finite group, the
homotopy is therefore given by (the links of) the intersecting
cells of codimension two (namely, the pentagon and commutativity
relations) together with the relation of serially and trivially crossing the same face
(namely, involutivity). $\Box$

\end{document}